\documentclass[12pt]{amsart}
\usepackage{amsmath}
\usepackage{amsfonts} 
\usepackage{amssymb}
\usepackage[all]{xy}
\usepackage{bbm}
\usepackage{epsfig}
\usepackage{url}
\usepackage{color}
\usepackage{epic}
\usepackage{eepic}
\usepackage{graphicx}

\newtheorem{theorem}{Theorem}
\newtheorem*{oniontheorem}{Onion--Skin Theorem}
\newtheorem{proposition}[theorem]{Proposition}
\newtheorem{exercise}[theorem]{Exercise}
\newtheorem{lemma}[theorem]{Lemma}
\theoremstyle{definition}

\newcommand{\R}{{\mathbb R}}
\newcommand{\Z}{{\mathbb Z}}
\renewcommand{\P}{{\mathbb P}}
\newcommand{\conv}{\mbox{\rm conv}}
\renewcommand{\int}{\operatorname{int}}
\newcommand{\suchthat}{\ : \ }
\newcommand{\slz}{\operatorname{SL_2 \Z}}
\newcommand{\pt}[1]{\mathbf{#1}}
\newcommand{\nv}{\pt{0}}
\newcommand{\segment}[2]{[
  {\pt{#1},\pt{#2}}]}
\newcommand{\sign}{\operatorname{sign}}
\newcommand{\transpose}{^t}

\begin{document}
\title{Lattice polygons and the number $\mathbf{2i+7}$}
\author{Christian Haase}
\address{Fachbereich Mathematik \& Informatik \\
  Freie Universit\"at Berlin \\
  14195 Berlin\\
  Germany}
\email{christian.haase@math.fu-berlin.de}
\urladdr{http://ehrhart.math.fu-berlin.de}
\thanks{Work by the first author was supported by NSF grant
  DMS--0200740 and by Emmy Noether fellowship HA 4383/1 of the German
  Research Foundation (DFG).
  Work by the second author was supported by the FWF in the frame
  of the SFB F013 and the project P15551.}
\author{Josef Schicho}
\address{RICAM \\
  Austrian Academy of Sciences \\
  Altenberger Stra\ss e 69 \\
  4040 Linz \\
  Austria}
\email{Josef.Schicho@oeaw.ac.at}
\urladdr{http://www.ricam.oeaw.ac.at/research/symcomp}

\begin{abstract}
  In this note we classify all triples $(a,b,i)$ such that there is a
  convex lattice polygon $P$ with area $a$ which has $b$ and $i$
  lattice points on the boundary and in the interior, respectively.
  The crucial lemma for the classification is the necessity of $b \le
  2 \ i + 7$. We sketch three proofs of this fact: the original one by
  Scott~\cite{Scott}, an elementary one, and one using algebraic
  geometry.

  As a refinement, we introduce an onion skin parameter $\ell$: how
  many nested polygons does $P$ contain? and give sharper bounds.
\end{abstract}

\maketitle
\setcounter{section}{-1}
\section{Introduction}\noindent
\subsection{How it all began}
When the second author translated a result on algebraic surfaces into
the language of lattice polygons using toric geometry, he obtained an
inequality for lattice polygons. This inequality had originally been
discovered by Scott~\cite{Scott}. The first author then found a third
proof. 
Subsequently, both authors
went through a phase of polygon addiction. Once you get started to
draw lattice polygons on graph paper and to discover relations between
their numerical invariants, it is not so easy to stop! (The gentle
reader has been warned.)

Thus, it was just unavoidable that the authors came up with new
inequalities: Scott's inequality can be sharpened if one takes another
invariant into account, which is defined by peeling off the skins of
the polygons like an onion (see Section~\ref{sec:onions}).

\subsection{Lattice polygons}
We want to study convex lattice polygons: convex polygons all whose
vertices have integral coordinates. As it turns out, we need to
consider nonconvex polygons as well. Even nonsimple polygons
-- polygons with self intersection -- will prove useful later on.
In what follows, we will abbreviate
\begin{center}
  ``polygon'' := ``convex lattice polygon'',
\end{center}
and we will emphasize when we allow nonintegral or nonsimple
situations. 
\begin{figure}[htbp]
  \begin{center}
\setlength{\unitlength}{0.00041667in}
\begingroup\makeatletter\ifx\SetFigFont\undefined%
\gdef\SetFigFont#1#2#3#4#5{%
  \reset@font\fontsize{#1}{#2pt}%
  \fontfamily{#3}\fontseries{#4}\fontshape{#5}%
  \selectfont}%
\fi\endgroup%
{\renewcommand{\dashlinestretch}{30}
\begin{picture}(7366,2581)(0,-10)
\put(7283,1283){\blacken\ellipse{150}{150}}
\put(7283,1283){\ellipse{150}{150}}
\put(7283,1883){\blacken\ellipse{150}{150}}
\put(7283,1883){\ellipse{150}{150}}
\put(7283,683){\blacken\ellipse{150}{150}}
\put(7283,683){\ellipse{150}{150}}
\put(6083,1283){\blacken\ellipse{150}{150}}
\put(6083,1283){\ellipse{150}{150}}
\put(6083,1883){\blacken\ellipse{150}{150}}
\put(6083,1883){\ellipse{150}{150}}
\put(6083,683){\blacken\ellipse{150}{150}}
\put(6083,683){\ellipse{150}{150}}
\put(6683,1283){\blacken\ellipse{150}{150}}
\put(6683,1283){\ellipse{150}{150}}
\put(6683,1883){\blacken\ellipse{150}{150}}
\put(6683,1883){\ellipse{150}{150}}
\put(6683,683){\blacken\ellipse{150}{150}}
\put(6683,683){\ellipse{150}{150}}
\thicklines
\path(7283,1283)(6083,1883)
\path(6324.495,1829.334)(6083.000,1883.000)(6270.830,1722.003)
\path(6083,1883)(6683,1883)
\path(6443.000,1823.000)(6683.000,1883.000)(6443.000,1943.000)
\path(6683,1883)(6083,683)
\path(6136.666,924.495)(6083.000,683.000)(6243.997,870.830)
\path(6083,683)(7283,1283)
\path(7095.170,1122.003)(7283.000,1283.000)(7041.505,1229.334)
\thinlines
\put(3083,1283){\blacken\ellipse{150}{150}}
\put(3083,1283){\ellipse{150}{150}}
\put(3083,1883){\blacken\ellipse{150}{150}}
\put(3083,1883){\ellipse{150}{150}}
\put(3083,683){\blacken\ellipse{150}{150}}
\put(3083,683){\ellipse{150}{150}}
\put(3683,1283){\blacken\ellipse{150}{150}}
\put(3683,1283){\ellipse{150}{150}}
\put(3683,1883){\blacken\ellipse{150}{150}}
\put(3683,1883){\ellipse{150}{150}}
\put(3683,683){\blacken\ellipse{150}{150}}
\put(3683,683){\ellipse{150}{150}}
\put(3683,2483){\blacken\ellipse{150}{150}}
\put(3683,2483){\ellipse{150}{150}}
\put(3083,2483){\blacken\ellipse{150}{150}}
\put(3083,2483){\ellipse{150}{150}}
\put(3683,83){\blacken\ellipse{150}{150}}
\put(3683,83){\ellipse{150}{150}}
\put(3083,83){\blacken\ellipse{150}{150}}
\put(3083,83){\ellipse{150}{150}}
\put(83,1283){\blacken\ellipse{150}{150}}
\put(83,1283){\ellipse{150}{150}}
\put(83,1883){\blacken\ellipse{150}{150}}
\put(83,1883){\ellipse{150}{150}}
\put(83,683){\blacken\ellipse{150}{150}}
\put(83,683){\ellipse{150}{150}}
\put(683,1283){\blacken\ellipse{150}{150}}
\put(683,1283){\ellipse{150}{150}}
\put(683,1883){\blacken\ellipse{150}{150}}
\put(683,1883){\ellipse{150}{150}}
\put(683,683){\blacken\ellipse{150}{150}}
\put(683,683){\ellipse{150}{150}}
\put(683,2483){\blacken\ellipse{150}{150}}
\put(683,2483){\ellipse{150}{150}}
\put(83,2483){\blacken\ellipse{150}{150}}
\put(83,2483){\ellipse{150}{150}}
\put(683,83){\blacken\ellipse{150}{150}}
\put(683,83){\ellipse{150}{150}}
\put(83,83){\blacken\ellipse{150}{150}}
\put(83,83){\ellipse{150}{150}}
\put(1283,1283){\blacken\ellipse{150}{150}}
\put(1283,1283){\ellipse{150}{150}}
\put(1283,1883){\blacken\ellipse{150}{150}}
\put(1283,1883){\ellipse{150}{150}}
\put(1283,683){\blacken\ellipse{150}{150}}
\put(1283,683){\ellipse{150}{150}}
\put(1883,1283){\blacken\ellipse{150}{150}}
\put(1883,1283){\ellipse{150}{150}}
\put(1883,1883){\blacken\ellipse{150}{150}}
\put(1883,1883){\ellipse{150}{150}}
\put(1883,683){\blacken\ellipse{150}{150}}
\put(1883,683){\ellipse{150}{150}}
\put(1883,2483){\blacken\ellipse{150}{150}}
\put(1883,2483){\ellipse{150}{150}}
\put(1283,2483){\blacken\ellipse{150}{150}}
\put(1283,2483){\ellipse{150}{150}}
\put(1883,83){\blacken\ellipse{150}{150}}
\put(1883,83){\ellipse{150}{150}}
\put(1283,83){\blacken\ellipse{150}{150}}
\put(1283,83){\ellipse{150}{150}}
\put(4283,1283){\blacken\ellipse{150}{150}}
\put(4283,1283){\ellipse{150}{150}}
\put(4283,1883){\blacken\ellipse{150}{150}}
\put(4283,1883){\ellipse{150}{150}}
\put(4283,683){\blacken\ellipse{150}{150}}
\put(4283,683){\ellipse{150}{150}}
\put(4283,2483){\blacken\ellipse{150}{150}}
\put(4283,2483){\ellipse{150}{150}}
\put(4883,2483){\blacken\ellipse{150}{150}}
\put(4883,2483){\ellipse{150}{150}}
\put(4883,1883){\blacken\ellipse{150}{150}}
\put(4883,1883){\ellipse{150}{150}}
\put(4883,1283){\blacken\ellipse{150}{150}}
\put(4883,1283){\ellipse{150}{150}}
\put(4883,683){\blacken\ellipse{150}{150}}
\put(4883,683){\ellipse{150}{150}}
\put(4883,83){\blacken\ellipse{150}{150}}
\put(4883,83){\ellipse{150}{150}}
\put(4283,83){\blacken\ellipse{150}{150}}
\put(4283,83){\ellipse{150}{150}}
\thicklines
\path(4883,2483)(4283,1883)
\path(4410.279,2095.132)(4283.000,1883.000)(4495.132,2010.279)
\path(4283,1883)(3083,1883)
\path(3323.000,1943.000)(3083.000,1883.000)(3323.000,1823.000)
\path(3083,1883)(4883,83)
\path(4670.868,210.279)(4883.000,83.000)(4755.721,295.132)
\path(4883,83)(4883,2483)
\path(4943.000,2243.000)(4883.000,2483.000)(4823.000,2243.000)
\path(1583,83)(83,83)
\path(323.000,143.000)(83.000,83.000)(323.000,23.000)
\path(383,2258)(1883,1283)
\path(1649.074,1363.491)(1883.000,1283.000)(1714.473,1464.104)
\path(83,83)(373,2187)
\path(399.668,1941.055)(373.000,2187.000)(280.792,1957.440)
\path(1883,1283)(1592,124)
\path(1592.251,371.386)(1592.000,124.000)(1708.639,342.164)
\end{picture}
}
    \caption{Polygons -- convex, lattice, nonsimple.}
    \label{fig:polygons}
  \end{center}
\end{figure}
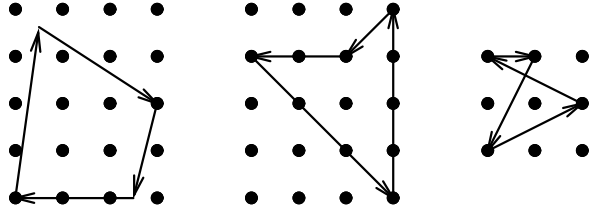

Denote the area enclosed by a polygon $P$ by $a=a(P)$, 
the number of lattice points on the boundary by 
$b=b(P)$, and the number of lattice points strictly inside of $P$ by
$i=i(P)$. 
A classic result relates these data.
\begin{theorem}[Pick's Formula~\cite{pick}]
  \begin{equation} \label{eq:pick}
    a = i + \frac{b}{2} - 1
  \end{equation}
\end{theorem}
A thorough discussion of this theorem -- including an
application in forest industry! -- can be found in
\cite{Gruenbaum_Shepard:93}.
Pick's theorem is not the only relation between the three parameters
$a$, $b$, $i$ of 
polygons.
There is the rather obvious constraint $b\ge 3$.
From Pick's formula we obtain immediately $a \ge i +
\frac{1}{2}$ and $a \ge \frac{b}{2} - 1$.
Are there other constraints?
For the sake of suspense, we do not want to reveal the final
inequalities just yet. We refer the impatient reader to the 
conclusion in Section~\ref{sec:conclusion} which contains a 
summary of the main results.
\subsection{Lattice equivalence} \label{sec:equivalence}
Clearly, the area $a(P)$ is invariant under rigid motions of the
plane. On the other hand, the numbers $i(P)$, $b(P)$ are not concepts
of Euclidean geometry, because they are not preserved by rigid
motions. But they {\em are} preserved under {\em
  lattice equivalences\/}: affine maps $\Phi \colon \R^2 \rightarrow
\R^2$ of the plane that restrict to isomorphisms of the lattice
$\Z^2$.

Orientation preserving lattice equivalences form a group, the semi
direct product $\slz \ltimes \Z^2$.

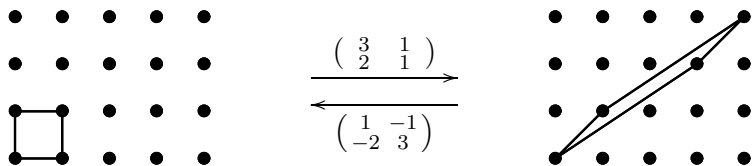
\begin{figure}[htbp]
  \centering
\setlength{\unitlength}{0.00041667in}
\begingroup\makeatletter\ifx\SetFigFont\undefined%
\gdef\SetFigFont#1#2#3#4#5{%
  \reset@font\fontsize{#1}{#2pt}%
  \fontfamily{#3}\fontseries{#4}\fontshape{#5}%
  \selectfont}%
\fi\endgroup%
{\renewcommand{\dashlinestretch}{30}
\begin{picture}(2566,1981)(0,-10)
\thicklines
\whiten\path(683,83)(83,83)(83,683)
	(683,683)(683,83)
\path(683,83)(83,83)(83,683)
	(683,683)(683,83)
\thinlines
\put(1883,1883){\blacken\ellipse{150}{150}}
\put(1883,1883){\ellipse{150}{150}}
\put(1283,1883){\blacken\ellipse{150}{150}}
\put(1283,1883){\ellipse{150}{150}}
\put(1883,683){\blacken\ellipse{150}{150}}
\put(1883,683){\ellipse{150}{150}}
\put(683,683){\blacken\ellipse{150}{150}}
\put(683,683){\ellipse{150}{150}}
\put(1883,1283){\blacken\ellipse{150}{150}}
\put(1883,1283){\ellipse{150}{150}}
\put(2483,1283){\blacken\ellipse{150}{150}}
\put(2483,1283){\ellipse{150}{150}}
\put(1283,1283){\blacken\ellipse{150}{150}}
\put(1283,1283){\ellipse{150}{150}}
\put(83,683){\blacken\ellipse{150}{150}}
\put(83,683){\ellipse{150}{150}}
\put(683,83){\blacken\ellipse{150}{150}}
\put(683,83){\ellipse{150}{150}}
\put(1283,83){\blacken\ellipse{150}{150}}
\put(1283,83){\ellipse{150}{150}}
\put(683,1283){\blacken\ellipse{150}{150}}
\put(683,1283){\ellipse{150}{150}}
\put(83,83){\blacken\ellipse{150}{150}}
\put(83,83){\ellipse{150}{150}}
\put(2483,1883){\blacken\ellipse{150}{150}}
\put(2483,1883){\ellipse{150}{150}}
\put(1883,83){\blacken\ellipse{150}{150}}
\put(1883,83){\ellipse{150}{150}}
\put(2483,83){\blacken\ellipse{150}{150}}
\put(2483,83){\ellipse{150}{150}}
\put(2483,683){\blacken\ellipse{150}{150}}
\put(2483,683){\ellipse{150}{150}}
\put(83,1283){\blacken\ellipse{150}{150}}
\put(83,1283){\ellipse{150}{150}}
\put(83,1883){\blacken\ellipse{150}{150}}
\put(83,1883){\ellipse{150}{150}}
\put(683,1883){\blacken\ellipse{150}{150}}
\put(683,1883){\ellipse{150}{150}}
\put(1283,683){\blacken\ellipse{150}{150}}
\put(1283,683){\ellipse{150}{150}}
\end{picture}
}
  \begin{minipage}{1.5in}
      \qquad \raisebox{25mm}{$\xymatrix{
        \ar@<1ex>[rr]^{%
          \left(\begin{smallmatrix}
            \ 3\ &\ 1\ \\2&1
          \end{smallmatrix}\right)
        }
        &&
        \ar@<1ex>[ll]^{%
          \left(\begin{smallmatrix}
            1&-1\\-2&3
          \end{smallmatrix}\right)
        }        
      }$}
  \end{minipage}
\setlength{\unitlength}{0.00041667in}
\begingroup\makeatletter\ifx\SetFigFont\undefined%
\gdef\SetFigFont#1#2#3#4#5{%
  \reset@font\fontsize{#1}{#2pt}%
  \fontfamily{#3}\fontseries{#4}\fontshape{#5}%
  \selectfont}%
\fi\endgroup%
{\renewcommand{\dashlinestretch}{30}
\begin{picture}(2566,1981)(0,-10)
\thicklines
\whiten\path(1883,1283)(83,83)(683,683)
	(2483,1883)(1883,1283)
\path(1883,1283)(83,83)(683,683)
	(2483,1883)(1883,1283)
\thinlines
\put(1883,1883){\blacken\ellipse{150}{150}}
\put(1883,1883){\ellipse{150}{150}}
\put(1283,1883){\blacken\ellipse{150}{150}}
\put(1283,1883){\ellipse{150}{150}}
\put(1883,683){\blacken\ellipse{150}{150}}
\put(1883,683){\ellipse{150}{150}}
\put(683,683){\blacken\ellipse{150}{150}}
\put(683,683){\ellipse{150}{150}}
\put(1283,683){\blacken\ellipse{150}{150}}
\put(1283,683){\ellipse{150}{150}}
\put(1883,1283){\blacken\ellipse{150}{150}}
\put(1883,1283){\ellipse{150}{150}}
\put(2483,1283){\blacken\ellipse{150}{150}}
\put(2483,1283){\ellipse{150}{150}}
\put(1283,1283){\blacken\ellipse{150}{150}}
\put(1283,1283){\ellipse{150}{150}}
\put(83,683){\blacken\ellipse{150}{150}}
\put(83,683){\ellipse{150}{150}}
\put(683,83){\blacken\ellipse{150}{150}}
\put(683,83){\ellipse{150}{150}}
\put(1283,83){\blacken\ellipse{150}{150}}
\put(1283,83){\ellipse{150}{150}}
\put(683,1283){\blacken\ellipse{150}{150}}
\put(683,1283){\ellipse{150}{150}}
\put(83,83){\blacken\ellipse{150}{150}}
\put(83,83){\ellipse{150}{150}}
\put(2483,1883){\blacken\ellipse{150}{150}}
\put(2483,1883){\ellipse{150}{150}}
\put(1883,83){\blacken\ellipse{150}{150}}
\put(1883,83){\ellipse{150}{150}}
\put(2483,83){\blacken\ellipse{150}{150}}
\put(2483,83){\ellipse{150}{150}}
\put(2483,683){\blacken\ellipse{150}{150}}
\put(2483,683){\ellipse{150}{150}}
\put(83,1283){\blacken\ellipse{150}{150}}
\put(83,1283){\ellipse{150}{150}}
\put(83,1883){\blacken\ellipse{150}{150}}
\put(83,1883){\ellipse{150}{150}}
\put(683,1883){\blacken\ellipse{150}{150}}
\put(683,1883){\ellipse{150}{150}}
\end{picture}
}
  \vspace{-15mm}
  \caption{Two lattice equivalent quadrangles.}
  \label{fig:square}
\end{figure}

So $\Phi$ has the form $\Phi(\pt{x}) = A \pt{x} + \pt{y}$ for a matrix
$A$, and a vector $\pt{y}$. The lattice preservation property
$\Phi(\Z^2)=\Z^2$ implies that both $A$ and $\pt{y}$ have
integral entries, and the same is true for the inverse transformation
$\Phi^{-1}(\pt{x}) = A^{-1} \pt{x} - A^{-1} \pt{y}$. Hence $\det A =
\pm 1$, and $a(P)$ is preserved under $\Phi$ as well. 

In all our arguments, we will treat lattice equivalent polygons as
indistinguishable. For example, the quadrangle in
Figure~\ref{fig:square} on the right looks to us like a perfect
square. We see that angles and Euclidian lengths are not preserved.
A lattice geometric substitute for the length of a lattice line
segment is the number of lattice points it contains minus one. 
In this sense, $b$ is the perimeter of $P$.
Here is an exercise that helps to get a feeling for what
lattice equivalences can do and cannot do.\footnote{Hint: if you know
  {\em Chinese\/}, not much {\em remains\/} to be done.} 
\begin{exercise} \label{ex:angle}
  Given a vertex $\pt{x}$ of a polygon $P$, show that there is a
  unique orientation preserving lattice equivalence $\Phi$ so that
  \begin{itemize}
  \item $\Phi(\pt{x})=(0,0)\transpose$, and
  \item 
    there are (necessarily unique) coprime $0 < p \le q$ so that
    the segments $[(1,0)\transpose,(0,0)\transpose]$ and
    $[(0,0)\transpose, (-p,q)\transpose]$ are contained in edges of
    $\Phi(P)$.
  \end{itemize}
\end{exercise}
\subsection{Why algebraic geometry?}
\begin{figure}[htbp]
  \centering
  \includegraphics[height=6cm,width=8cm]{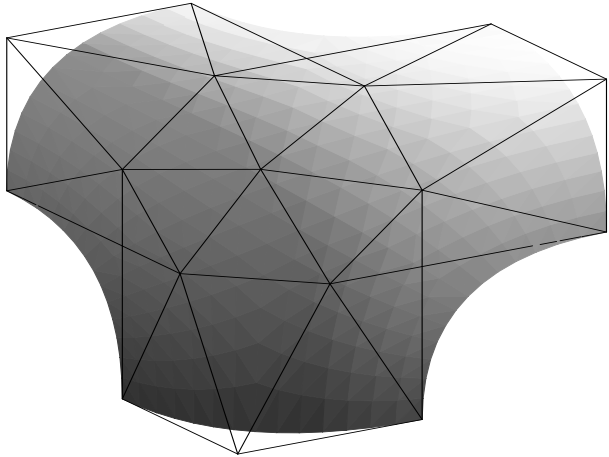}
  \caption{Toric surface with hexagonal control structure.}
  \label{fig:rimas}
\end{figure}

Toric geometry is a powerful link connecting discrete and algebraic
geometry (see e.g. \cite{Sturmfels:99}).
At the heart of this link is the simple correspondence
\begin{eqnarray*}
  \text{ lattice point \qquad } & & \text{ \qquad Laurent monomial } \\
  \pt{p}=(p_1,\ldots,p_m) \in \Z^m & \longleftrightarrow &
  \pt{x}^{\pt{p}} = x_1^{p_1} \cdot\ldots\cdot x_m^{p_m} \in
  \mathbb{C}[x_1^{\pm1},\ldots,x_m^{\pm1}]
\end{eqnarray*} 
It was invented by M.~Demazure~\cite{Demazure:70} for a totally
different purpose (to study algebraic subgroups of the Cremona group)
in algebraic geometry. R.~Stanley used it in combinatorics to classify
the possible face numbers of simplicial convex polytopes~\cite{gThmS}.
R.~Krasauskas~\cite{Krasauskas:02} used it in geometric modeling to
construct surfaces with new control structure (see
Figure~\ref{fig:rimas}). 

For any polygon $P$, the Laurent monomials corresponding to its
lattice points define a {\em toric surface\/} $X_P$ in a projective
space of dimension $b+i-1$ as follows.
Number the lattice points $P \cap \Z^2=\{\pt{p}_0,\dots,\pt{p}_n\}$
(where $n=b+i-1$). Then $X_P$ is the closure of the image of the map
$({\mathbb C}^*)^2\to{\mathbb P}^n$ defined by $\pt{x} \mapsto
(\pt{x}^{\pt{p}_0} : \dots : \pt{x}^{\pt{p}_n})$.
Lattice equivalent polygons define the same toric surface. 

As to be expected, there is a dictionary translating toric geometry
to lattice geometry:
the degree of the toric surface is equal to twice
the area, and the number of interior points is equal to
the {\em sectional genus\/} of the surface.
For instance, let $\Gamma$ be the triangle with corners
$(0,0)^t,(1,2)^t,(2,1)^t$ (it has one interior point $(1,1)^t$).
Then the toric surface is given by $(1 : x_1x_2^2 : x_1^2x_2 : x_1x_2)
\in \P^3$, its degree is 3 -- this is also reflected by its implicit
equation $y_1y_2y_3-y_4^3=0$, which has also degree~$3$ -- 
and its sectional genus is $1$, i.e., if we intersect with a generic
hyperplane in $\P^3$, we obtain a genus one Riemann surface.

In the context of toric geometry, Pick's formula appears 
as a consequence of the Riemann-Roch Theorem.
 
\section{Examples}
\begin{figure}[htbp]
  \centering
\setlength{\unitlength}{0.00041667in}
\begingroup\makeatletter\ifx\SetFigFont\undefined%
\gdef\SetFigFont#1#2#3#4#5{%
  \reset@font\fontsize{#1}{#2pt}%
  \fontfamily{#3}\fontseries{#4}\fontshape{#5}%
  \selectfont}%
\fi\endgroup%
{\renewcommand{\dashlinestretch}{30}
\begin{picture}(9766,3181)(0,-10)
\put(1883,983){\blacken\ellipse{150}{150}}
\put(1883,983){\ellipse{150}{150}}
\put(2483,983){\blacken\ellipse{150}{150}}
\put(2483,983){\ellipse{150}{150}}
\put(3083,983){\blacken\ellipse{150}{150}}
\put(3083,983){\ellipse{150}{150}}
\put(683,983){\blacken\ellipse{150}{150}}
\put(683,983){\ellipse{150}{150}}
\put(1283,983){\blacken\ellipse{150}{150}}
\put(1283,983){\ellipse{150}{150}}
\put(83,2183){\blacken\ellipse{150}{150}}
\put(83,2183){\ellipse{150}{150}}
\put(3683,2183){\blacken\ellipse{150}{150}}
\put(3683,2183){\ellipse{150}{150}}
\put(3683,1583){\blacken\ellipse{150}{150}}
\put(3683,1583){\ellipse{150}{150}}
\put(83,983){\blacken\ellipse{150}{150}}
\put(83,983){\ellipse{150}{150}}
\put(83,1583){\blacken\ellipse{150}{150}}
\put(83,1583){\ellipse{150}{150}}
\put(1883,1583){\blacken\ellipse{150}{150}}
\put(1883,1583){\ellipse{150}{150}}
\put(2483,1583){\blacken\ellipse{150}{150}}
\put(2483,1583){\ellipse{150}{150}}
\put(3083,1583){\blacken\ellipse{150}{150}}
\put(3083,1583){\ellipse{150}{150}}
\put(683,1583){\blacken\ellipse{150}{150}}
\put(683,1583){\ellipse{150}{150}}
\put(1283,1583){\blacken\ellipse{150}{150}}
\put(1283,1583){\ellipse{150}{150}}
\put(1883,2183){\blacken\ellipse{150}{150}}
\put(1883,2183){\ellipse{150}{150}}
\put(2483,2183){\blacken\ellipse{150}{150}}
\put(2483,2183){\ellipse{150}{150}}
\put(3083,2183){\blacken\ellipse{150}{150}}
\put(3083,2183){\ellipse{150}{150}}
\put(683,2183){\blacken\ellipse{150}{150}}
\put(683,2183){\ellipse{150}{150}}
\put(1283,2183){\blacken\ellipse{150}{150}}
\put(1283,2183){\ellipse{150}{150}}
\put(3683,983){\blacken\ellipse{150}{150}}
\put(3683,983){\ellipse{150}{150}}
\put(8483,1883){\blacken\ellipse{150}{150}}
\put(8483,1883){\ellipse{150}{150}}
\put(9083,1883){\blacken\ellipse{150}{150}}
\put(9083,1883){\ellipse{150}{150}}
\put(9683,1883){\blacken\ellipse{150}{150}}
\put(9683,1883){\ellipse{150}{150}}
\put(7283,1883){\blacken\ellipse{150}{150}}
\put(7283,1883){\ellipse{150}{150}}
\put(7883,1883){\blacken\ellipse{150}{150}}
\put(7883,1883){\ellipse{150}{150}}
\put(6683,1883){\blacken\ellipse{150}{150}}
\put(6683,1883){\ellipse{150}{150}}
\put(8483,2483){\blacken\ellipse{150}{150}}
\put(8483,2483){\ellipse{150}{150}}
\put(9083,2483){\blacken\ellipse{150}{150}}
\put(9083,2483){\ellipse{150}{150}}
\put(9683,2483){\blacken\ellipse{150}{150}}
\put(9683,2483){\ellipse{150}{150}}
\put(7283,2483){\blacken\ellipse{150}{150}}
\put(7283,2483){\ellipse{150}{150}}
\put(7883,2483){\blacken\ellipse{150}{150}}
\put(7883,2483){\ellipse{150}{150}}
\put(6683,2483){\blacken\ellipse{150}{150}}
\put(6683,2483){\ellipse{150}{150}}
\put(8483,3083){\blacken\ellipse{150}{150}}
\put(8483,3083){\ellipse{150}{150}}
\put(9083,3083){\blacken\ellipse{150}{150}}
\put(9083,3083){\ellipse{150}{150}}
\put(9683,3083){\blacken\ellipse{150}{150}}
\put(9683,3083){\ellipse{150}{150}}
\put(7283,3083){\blacken\ellipse{150}{150}}
\put(7283,3083){\ellipse{150}{150}}
\put(7883,3083){\blacken\ellipse{150}{150}}
\put(7883,3083){\ellipse{150}{150}}
\put(6683,3083){\blacken\ellipse{150}{150}}
\put(6683,3083){\ellipse{150}{150}}
\put(8483,1283){\blacken\ellipse{150}{150}}
\put(8483,1283){\ellipse{150}{150}}
\put(9083,1283){\blacken\ellipse{150}{150}}
\put(9083,1283){\ellipse{150}{150}}
\put(9683,1283){\blacken\ellipse{150}{150}}
\put(9683,1283){\ellipse{150}{150}}
\put(8483,83){\blacken\ellipse{150}{150}}
\put(8483,83){\ellipse{150}{150}}
\put(9083,83){\blacken\ellipse{150}{150}}
\put(9083,83){\ellipse{150}{150}}
\put(9683,83){\blacken\ellipse{150}{150}}
\put(9683,83){\ellipse{150}{150}}
\put(8483,683){\blacken\ellipse{150}{150}}
\put(8483,683){\ellipse{150}{150}}
\put(9083,683){\blacken\ellipse{150}{150}}
\put(9083,683){\ellipse{150}{150}}
\put(9683,683){\blacken\ellipse{150}{150}}
\put(9683,683){\ellipse{150}{150}}
\put(7283,1283){\blacken\ellipse{150}{150}}
\put(7283,1283){\ellipse{150}{150}}
\put(7283,83){\blacken\ellipse{150}{150}}
\put(7283,83){\ellipse{150}{150}}
\put(7283,683){\blacken\ellipse{150}{150}}
\put(7283,683){\ellipse{150}{150}}
\put(7883,1283){\blacken\ellipse{150}{150}}
\put(7883,1283){\ellipse{150}{150}}
\put(7883,683){\blacken\ellipse{150}{150}}
\put(7883,683){\ellipse{150}{150}}
\put(6683,1283){\blacken\ellipse{150}{150}}
\put(6683,1283){\ellipse{150}{150}}
\put(6683,83){\blacken\ellipse{150}{150}}
\put(6683,83){\ellipse{150}{150}}
\put(6683,683){\blacken\ellipse{150}{150}}
\put(6683,683){\ellipse{150}{150}}
\put(7883,83){\blacken\ellipse{150}{150}}
\put(7883,83){\ellipse{150}{150}}
\thicklines
\path(83,983)(683,2183)(3683,1583)(83,983)
\path(7283,3083)(9683,83)(6683,683)(7283,3083)
\end{picture}
}
  \caption{$b=3$ and $a \gg 0$}
  \label{fig:bequal3}
\end{figure}
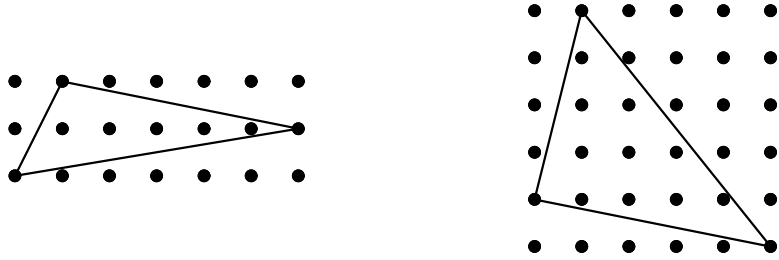

\noindent
%
Let us approach the question which parameters are possible for
polygons by looking at some examples. Can we bound $i$ or $a$
in terms of $b$\,?   Figure~\ref{fig:bequal3} shows examples with $b=3$
and arbitrarily high $a$ and $i$. So there is no lattice geometric
analogue of the isoperimetric inequality.

\begin{figure}[htbp]
  \centering
\setlength{\unitlength}{0.00041667in}
\begingroup\makeatletter\ifx\SetFigFont\undefined%
\gdef\SetFigFont#1#2#3#4#5{%
  \reset@font\fontsize{#1}{#2pt}%
  \fontfamily{#3}\fontseries{#4}\fontshape{#5}%
  \selectfont}%
\fi\endgroup%
{\renewcommand{\dashlinestretch}{30}
\begin{picture}(3766,781)(0,-10)
\put(83,683){\blacken\ellipse{150}{150}}
\put(83,683){\ellipse{150}{150}}
\put(3683,683){\blacken\ellipse{150}{150}}
\put(3683,683){\ellipse{150}{150}}
\put(3683,83){\blacken\ellipse{150}{150}}
\put(3683,83){\ellipse{150}{150}}
\put(83,83){\blacken\ellipse{150}{150}}
\put(83,83){\ellipse{150}{150}}
\put(1883,83){\blacken\ellipse{150}{150}}
\put(1883,83){\ellipse{150}{150}}
\put(2483,83){\blacken\ellipse{150}{150}}
\put(2483,83){\ellipse{150}{150}}
\put(3083,83){\blacken\ellipse{150}{150}}
\put(3083,83){\ellipse{150}{150}}
\put(683,83){\blacken\ellipse{150}{150}}
\put(683,83){\ellipse{150}{150}}
\put(1283,83){\blacken\ellipse{150}{150}}
\put(1283,83){\ellipse{150}{150}}
\put(1883,683){\blacken\ellipse{150}{150}}
\put(1883,683){\ellipse{150}{150}}
\put(2483,683){\blacken\ellipse{150}{150}}
\put(2483,683){\ellipse{150}{150}}
\put(3083,683){\blacken\ellipse{150}{150}}
\put(3083,683){\ellipse{150}{150}}
\put(683,683){\blacken\ellipse{150}{150}}
\put(683,683){\ellipse{150}{150}}
\put(1283,683){\blacken\ellipse{150}{150}}
\put(1283,683){\ellipse{150}{150}}
\thicklines
\path(83,683)(83,83)(3683,83)(83,683)
\end{picture}
}
  \caption{$i=0$ and $b \gg 0$.}
  \label{fig:iequal0}
\end{figure}
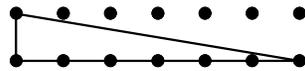

What about bounds in the opposite direction? Can we bound $b$ in terms
of $i$\,? Well, there is the family of Figure~\ref{fig:iequal0} with
$i=0$ and arbitrary $b$.

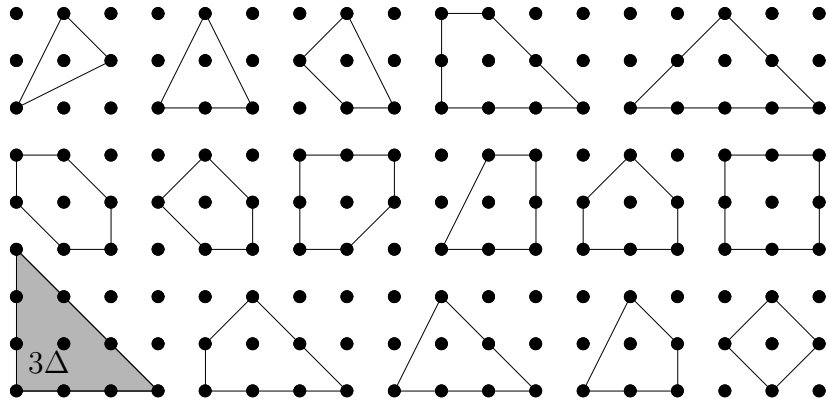
\begin{figure}[htbp]
  \centering
\setlength{\unitlength}{0.00041667in}
\begingroup\makeatletter\ifx\SetFigFont\undefined%
\gdef\SetFigFont#1#2#3#4#5{%
  \reset@font\fontsize{#1}{#2pt}%
  \fontfamily{#3}\fontseries{#4}\fontshape{#5}%
  \selectfont}%
\fi\endgroup%
{\renewcommand{\dashlinestretch}{30}
\begin{picture}(10366,4981)(0,-10)
\texture{88555555 55000000 555555 55000000 555555 55000000 555555 55000000 
	555555 55000000 555555 55000000 555555 55000000 555555 55000000 
	555555 55000000 555555 55000000 555555 55000000 555555 55000000 
	555555 55000000 555555 55000000 555555 55000000 555555 55000000 }
\shade\path(83,1883)(83,83)(1883,83)(83,1883)
\path(83,1883)(83,83)(1883,83)(83,1883)
\put(83,4883){\blacken\ellipse{150}{150}}
\put(83,4883){\ellipse{150}{150}}
\put(1283,4883){\blacken\ellipse{150}{150}}
\put(1283,4883){\ellipse{150}{150}}
\put(2483,4883){\blacken\ellipse{150}{150}}
\put(2483,4883){\ellipse{150}{150}}
\put(3683,4883){\blacken\ellipse{150}{150}}
\put(3683,4883){\ellipse{150}{150}}
\put(4883,4883){\blacken\ellipse{150}{150}}
\put(4883,4883){\ellipse{150}{150}}
\put(6083,4883){\blacken\ellipse{150}{150}}
\put(6083,4883){\ellipse{150}{150}}
\put(7283,4883){\blacken\ellipse{150}{150}}
\put(7283,4883){\ellipse{150}{150}}
\put(8483,4883){\blacken\ellipse{150}{150}}
\put(8483,4883){\ellipse{150}{150}}
\put(9683,4883){\blacken\ellipse{150}{150}}
\put(9683,4883){\ellipse{150}{150}}
\put(683,4883){\blacken\ellipse{150}{150}}
\put(683,4883){\ellipse{150}{150}}
\put(1883,4883){\blacken\ellipse{150}{150}}
\put(1883,4883){\ellipse{150}{150}}
\put(3083,4883){\blacken\ellipse{150}{150}}
\put(3083,4883){\ellipse{150}{150}}
\put(4283,4883){\blacken\ellipse{150}{150}}
\put(4283,4883){\ellipse{150}{150}}
\put(5483,4883){\blacken\ellipse{150}{150}}
\put(5483,4883){\ellipse{150}{150}}
\put(6683,4883){\blacken\ellipse{150}{150}}
\put(6683,4883){\ellipse{150}{150}}
\put(7883,4883){\blacken\ellipse{150}{150}}
\put(7883,4883){\ellipse{150}{150}}
\put(9083,4883){\blacken\ellipse{150}{150}}
\put(9083,4883){\ellipse{150}{150}}
\put(10283,4883){\blacken\ellipse{150}{150}}
\put(10283,4883){\ellipse{150}{150}}
\put(83,4283){\blacken\ellipse{150}{150}}
\put(83,4283){\ellipse{150}{150}}
\put(1283,4283){\blacken\ellipse{150}{150}}
\put(1283,4283){\ellipse{150}{150}}
\put(2483,4283){\blacken\ellipse{150}{150}}
\put(2483,4283){\ellipse{150}{150}}
\put(3683,4283){\blacken\ellipse{150}{150}}
\put(3683,4283){\ellipse{150}{150}}
\put(4883,4283){\blacken\ellipse{150}{150}}
\put(4883,4283){\ellipse{150}{150}}
\put(6083,4283){\blacken\ellipse{150}{150}}
\put(6083,4283){\ellipse{150}{150}}
\put(7283,4283){\blacken\ellipse{150}{150}}
\put(7283,4283){\ellipse{150}{150}}
\put(8483,4283){\blacken\ellipse{150}{150}}
\put(8483,4283){\ellipse{150}{150}}
\put(9683,4283){\blacken\ellipse{150}{150}}
\put(9683,4283){\ellipse{150}{150}}
\put(683,4283){\blacken\ellipse{150}{150}}
\put(683,4283){\ellipse{150}{150}}
\put(1883,4283){\blacken\ellipse{150}{150}}
\put(1883,4283){\ellipse{150}{150}}
\put(3083,4283){\blacken\ellipse{150}{150}}
\put(3083,4283){\ellipse{150}{150}}
\put(4283,4283){\blacken\ellipse{150}{150}}
\put(4283,4283){\ellipse{150}{150}}
\put(5483,4283){\blacken\ellipse{150}{150}}
\put(5483,4283){\ellipse{150}{150}}
\put(6683,4283){\blacken\ellipse{150}{150}}
\put(6683,4283){\ellipse{150}{150}}
\put(7883,4283){\blacken\ellipse{150}{150}}
\put(7883,4283){\ellipse{150}{150}}
\put(9083,4283){\blacken\ellipse{150}{150}}
\put(9083,4283){\ellipse{150}{150}}
\put(10283,4283){\blacken\ellipse{150}{150}}
\put(10283,4283){\ellipse{150}{150}}
\put(83,3683){\blacken\ellipse{150}{150}}
\put(83,3683){\ellipse{150}{150}}
\put(1283,3683){\blacken\ellipse{150}{150}}
\put(1283,3683){\ellipse{150}{150}}
\put(2483,3683){\blacken\ellipse{150}{150}}
\put(2483,3683){\ellipse{150}{150}}
\put(3683,3683){\blacken\ellipse{150}{150}}
\put(3683,3683){\ellipse{150}{150}}
\put(4883,3683){\blacken\ellipse{150}{150}}
\put(4883,3683){\ellipse{150}{150}}
\put(6083,3683){\blacken\ellipse{150}{150}}
\put(6083,3683){\ellipse{150}{150}}
\put(7283,3683){\blacken\ellipse{150}{150}}
\put(7283,3683){\ellipse{150}{150}}
\put(8483,3683){\blacken\ellipse{150}{150}}
\put(8483,3683){\ellipse{150}{150}}
\put(9683,3683){\blacken\ellipse{150}{150}}
\put(9683,3683){\ellipse{150}{150}}
\put(683,3683){\blacken\ellipse{150}{150}}
\put(683,3683){\ellipse{150}{150}}
\put(1883,3683){\blacken\ellipse{150}{150}}
\put(1883,3683){\ellipse{150}{150}}
\put(3083,3683){\blacken\ellipse{150}{150}}
\put(3083,3683){\ellipse{150}{150}}
\put(4283,3683){\blacken\ellipse{150}{150}}
\put(4283,3683){\ellipse{150}{150}}
\put(5483,3683){\blacken\ellipse{150}{150}}
\put(5483,3683){\ellipse{150}{150}}
\put(6683,3683){\blacken\ellipse{150}{150}}
\put(6683,3683){\ellipse{150}{150}}
\put(7883,3683){\blacken\ellipse{150}{150}}
\put(7883,3683){\ellipse{150}{150}}
\put(9083,3683){\blacken\ellipse{150}{150}}
\put(9083,3683){\ellipse{150}{150}}
\put(10283,3683){\blacken\ellipse{150}{150}}
\put(10283,3683){\ellipse{150}{150}}
\put(83,3083){\blacken\ellipse{150}{150}}
\put(83,3083){\ellipse{150}{150}}
\put(1283,3083){\blacken\ellipse{150}{150}}
\put(1283,3083){\ellipse{150}{150}}
\put(2483,3083){\blacken\ellipse{150}{150}}
\put(2483,3083){\ellipse{150}{150}}
\put(3683,3083){\blacken\ellipse{150}{150}}
\put(3683,3083){\ellipse{150}{150}}
\put(4883,3083){\blacken\ellipse{150}{150}}
\put(4883,3083){\ellipse{150}{150}}
\put(6083,3083){\blacken\ellipse{150}{150}}
\put(6083,3083){\ellipse{150}{150}}
\put(7283,3083){\blacken\ellipse{150}{150}}
\put(7283,3083){\ellipse{150}{150}}
\put(8483,3083){\blacken\ellipse{150}{150}}
\put(8483,3083){\ellipse{150}{150}}
\put(9683,3083){\blacken\ellipse{150}{150}}
\put(9683,3083){\ellipse{150}{150}}
\put(683,3083){\blacken\ellipse{150}{150}}
\put(683,3083){\ellipse{150}{150}}
\put(1883,3083){\blacken\ellipse{150}{150}}
\put(1883,3083){\ellipse{150}{150}}
\put(3083,3083){\blacken\ellipse{150}{150}}
\put(3083,3083){\ellipse{150}{150}}
\put(4283,3083){\blacken\ellipse{150}{150}}
\put(4283,3083){\ellipse{150}{150}}
\put(5483,3083){\blacken\ellipse{150}{150}}
\put(5483,3083){\ellipse{150}{150}}
\put(6683,3083){\blacken\ellipse{150}{150}}
\put(6683,3083){\ellipse{150}{150}}
\put(7883,3083){\blacken\ellipse{150}{150}}
\put(7883,3083){\ellipse{150}{150}}
\put(9083,3083){\blacken\ellipse{150}{150}}
\put(9083,3083){\ellipse{150}{150}}
\put(10283,3083){\blacken\ellipse{150}{150}}
\put(10283,3083){\ellipse{150}{150}}
\put(83,2483){\blacken\ellipse{150}{150}}
\put(83,2483){\ellipse{150}{150}}
\put(1283,2483){\blacken\ellipse{150}{150}}
\put(1283,2483){\ellipse{150}{150}}
\put(2483,2483){\blacken\ellipse{150}{150}}
\put(2483,2483){\ellipse{150}{150}}
\put(3683,2483){\blacken\ellipse{150}{150}}
\put(3683,2483){\ellipse{150}{150}}
\put(4883,2483){\blacken\ellipse{150}{150}}
\put(4883,2483){\ellipse{150}{150}}
\put(6083,2483){\blacken\ellipse{150}{150}}
\put(6083,2483){\ellipse{150}{150}}
\put(7283,2483){\blacken\ellipse{150}{150}}
\put(7283,2483){\ellipse{150}{150}}
\put(8483,2483){\blacken\ellipse{150}{150}}
\put(8483,2483){\ellipse{150}{150}}
\put(9683,2483){\blacken\ellipse{150}{150}}
\put(9683,2483){\ellipse{150}{150}}
\put(683,2483){\blacken\ellipse{150}{150}}
\put(683,2483){\ellipse{150}{150}}
\put(1883,2483){\blacken\ellipse{150}{150}}
\put(1883,2483){\ellipse{150}{150}}
\put(3083,2483){\blacken\ellipse{150}{150}}
\put(3083,2483){\ellipse{150}{150}}
\put(4283,2483){\blacken\ellipse{150}{150}}
\put(4283,2483){\ellipse{150}{150}}
\put(5483,2483){\blacken\ellipse{150}{150}}
\put(5483,2483){\ellipse{150}{150}}
\put(6683,2483){\blacken\ellipse{150}{150}}
\put(6683,2483){\ellipse{150}{150}}
\put(7883,2483){\blacken\ellipse{150}{150}}
\put(7883,2483){\ellipse{150}{150}}
\put(9083,2483){\blacken\ellipse{150}{150}}
\put(9083,2483){\ellipse{150}{150}}
\put(10283,2483){\blacken\ellipse{150}{150}}
\put(10283,2483){\ellipse{150}{150}}
\put(83,1883){\blacken\ellipse{150}{150}}
\put(83,1883){\ellipse{150}{150}}
\put(1283,1883){\blacken\ellipse{150}{150}}
\put(1283,1883){\ellipse{150}{150}}
\put(2483,1883){\blacken\ellipse{150}{150}}
\put(2483,1883){\ellipse{150}{150}}
\put(3683,1883){\blacken\ellipse{150}{150}}
\put(3683,1883){\ellipse{150}{150}}
\put(4883,1883){\blacken\ellipse{150}{150}}
\put(4883,1883){\ellipse{150}{150}}
\put(6083,1883){\blacken\ellipse{150}{150}}
\put(6083,1883){\ellipse{150}{150}}
\put(7283,1883){\blacken\ellipse{150}{150}}
\put(7283,1883){\ellipse{150}{150}}
\put(8483,1883){\blacken\ellipse{150}{150}}
\put(8483,1883){\ellipse{150}{150}}
\put(9683,1883){\blacken\ellipse{150}{150}}
\put(9683,1883){\ellipse{150}{150}}
\put(683,1883){\blacken\ellipse{150}{150}}
\put(683,1883){\ellipse{150}{150}}
\put(1883,1883){\blacken\ellipse{150}{150}}
\put(1883,1883){\ellipse{150}{150}}
\put(3083,1883){\blacken\ellipse{150}{150}}
\put(3083,1883){\ellipse{150}{150}}
\put(4283,1883){\blacken\ellipse{150}{150}}
\put(4283,1883){\ellipse{150}{150}}
\put(5483,1883){\blacken\ellipse{150}{150}}
\put(5483,1883){\ellipse{150}{150}}
\put(6683,1883){\blacken\ellipse{150}{150}}
\put(6683,1883){\ellipse{150}{150}}
\put(7883,1883){\blacken\ellipse{150}{150}}
\put(7883,1883){\ellipse{150}{150}}
\put(9083,1883){\blacken\ellipse{150}{150}}
\put(9083,1883){\ellipse{150}{150}}
\put(10283,1883){\blacken\ellipse{150}{150}}
\put(10283,1883){\ellipse{150}{150}}
\put(83,1283){\blacken\ellipse{150}{150}}
\put(83,1283){\ellipse{150}{150}}
\put(1283,1283){\blacken\ellipse{150}{150}}
\put(1283,1283){\ellipse{150}{150}}
\put(2483,1283){\blacken\ellipse{150}{150}}
\put(2483,1283){\ellipse{150}{150}}
\put(3683,1283){\blacken\ellipse{150}{150}}
\put(3683,1283){\ellipse{150}{150}}
\put(4883,1283){\blacken\ellipse{150}{150}}
\put(4883,1283){\ellipse{150}{150}}
\put(6083,1283){\blacken\ellipse{150}{150}}
\put(6083,1283){\ellipse{150}{150}}
\put(7283,1283){\blacken\ellipse{150}{150}}
\put(7283,1283){\ellipse{150}{150}}
\put(8483,1283){\blacken\ellipse{150}{150}}
\put(8483,1283){\ellipse{150}{150}}
\put(9683,1283){\blacken\ellipse{150}{150}}
\put(9683,1283){\ellipse{150}{150}}
\put(683,1283){\blacken\ellipse{150}{150}}
\put(683,1283){\ellipse{150}{150}}
\put(1883,1283){\blacken\ellipse{150}{150}}
\put(1883,1283){\ellipse{150}{150}}
\put(3083,1283){\blacken\ellipse{150}{150}}
\put(3083,1283){\ellipse{150}{150}}
\put(4283,1283){\blacken\ellipse{150}{150}}
\put(4283,1283){\ellipse{150}{150}}
\put(5483,1283){\blacken\ellipse{150}{150}}
\put(5483,1283){\ellipse{150}{150}}
\put(6683,1283){\blacken\ellipse{150}{150}}
\put(6683,1283){\ellipse{150}{150}}
\put(7883,1283){\blacken\ellipse{150}{150}}
\put(7883,1283){\ellipse{150}{150}}
\put(9083,1283){\blacken\ellipse{150}{150}}
\put(9083,1283){\ellipse{150}{150}}
\put(10283,1283){\blacken\ellipse{150}{150}}
\put(10283,1283){\ellipse{150}{150}}
\put(83,683){\blacken\ellipse{150}{150}}
\put(83,683){\ellipse{150}{150}}
\put(1283,683){\blacken\ellipse{150}{150}}
\put(1283,683){\ellipse{150}{150}}
\put(2483,683){\blacken\ellipse{150}{150}}
\put(2483,683){\ellipse{150}{150}}
\put(3683,683){\blacken\ellipse{150}{150}}
\put(3683,683){\ellipse{150}{150}}
\put(4883,683){\blacken\ellipse{150}{150}}
\put(4883,683){\ellipse{150}{150}}
\put(6083,683){\blacken\ellipse{150}{150}}
\put(6083,683){\ellipse{150}{150}}
\put(7283,683){\blacken\ellipse{150}{150}}
\put(7283,683){\ellipse{150}{150}}
\put(8483,683){\blacken\ellipse{150}{150}}
\put(8483,683){\ellipse{150}{150}}
\put(9683,683){\blacken\ellipse{150}{150}}
\put(9683,683){\ellipse{150}{150}}
\put(683,683){\blacken\ellipse{150}{150}}
\put(683,683){\ellipse{150}{150}}
\put(1883,683){\blacken\ellipse{150}{150}}
\put(1883,683){\ellipse{150}{150}}
\put(3083,683){\blacken\ellipse{150}{150}}
\put(3083,683){\ellipse{150}{150}}
\put(4283,683){\blacken\ellipse{150}{150}}
\put(4283,683){\ellipse{150}{150}}
\put(5483,683){\blacken\ellipse{150}{150}}
\put(5483,683){\ellipse{150}{150}}
\put(6683,683){\blacken\ellipse{150}{150}}
\put(6683,683){\ellipse{150}{150}}
\put(7883,683){\blacken\ellipse{150}{150}}
\put(7883,683){\ellipse{150}{150}}
\put(9083,683){\blacken\ellipse{150}{150}}
\put(9083,683){\ellipse{150}{150}}
\put(10283,683){\blacken\ellipse{150}{150}}
\put(10283,683){\ellipse{150}{150}}
\put(83,83){\blacken\ellipse{150}{150}}
\put(83,83){\ellipse{150}{150}}
\put(1283,83){\blacken\ellipse{150}{150}}
\put(1283,83){\ellipse{150}{150}}
\put(2483,83){\blacken\ellipse{150}{150}}
\put(2483,83){\ellipse{150}{150}}
\put(3683,83){\blacken\ellipse{150}{150}}
\put(3683,83){\ellipse{150}{150}}
\put(4883,83){\blacken\ellipse{150}{150}}
\put(4883,83){\ellipse{150}{150}}
\put(6083,83){\blacken\ellipse{150}{150}}
\put(6083,83){\ellipse{150}{150}}
\put(7283,83){\blacken\ellipse{150}{150}}
\put(7283,83){\ellipse{150}{150}}
\put(8483,83){\blacken\ellipse{150}{150}}
\put(8483,83){\ellipse{150}{150}}
\put(9683,83){\blacken\ellipse{150}{150}}
\put(9683,83){\ellipse{150}{150}}
\put(683,83){\blacken\ellipse{150}{150}}
\put(683,83){\ellipse{150}{150}}
\put(1883,83){\blacken\ellipse{150}{150}}
\put(1883,83){\ellipse{150}{150}}
\put(3083,83){\blacken\ellipse{150}{150}}
\put(3083,83){\ellipse{150}{150}}
\put(4283,83){\blacken\ellipse{150}{150}}
\put(4283,83){\ellipse{150}{150}}
\put(5483,83){\blacken\ellipse{150}{150}}
\put(5483,83){\ellipse{150}{150}}
\put(6683,83){\blacken\ellipse{150}{150}}
\put(6683,83){\ellipse{150}{150}}
\put(7883,83){\blacken\ellipse{150}{150}}
\put(7883,83){\ellipse{150}{150}}
\put(9083,83){\blacken\ellipse{150}{150}}
\put(9083,83){\ellipse{150}{150}}
\put(10283,83){\blacken\ellipse{150}{150}}
\put(10283,83){\ellipse{150}{150}}
\path(683,4883)(1283,4283)(83,3683)(683,4883)
\path(1883,3683)(2483,4883)(3083,3683)(1883,3683)
\path(3683,4283)(4283,3683)(4883,3683)
	(4283,4883)(3683,4283)
\path(5483,4883)(5483,3683)(7283,3683)
	(6083,4883)(5483,4883)
\path(7883,3683)(9083,4883)(10283,3683)(7883,3683)
\path(83,3083)(83,2483)(683,1883)
	(1283,1883)(1283,2483)(683,3083)(83,3083)
\path(1883,2483)(2483,3083)(3083,2483)
	(3083,1883)(2483,1883)(1883,2483)
\path(3683,1883)(3683,3083)(4883,3083)
	(4883,2483)(4283,1883)(3683,1883)
\path(5483,1883)(6083,3083)(6683,3083)
	(6683,1883)(5483,1883)
\path(7283,1883)(7283,2483)(7883,3083)
	(8483,2483)(8483,1883)(7283,1883)
\path(9083,1883)(9083,3083)(10283,3083)
	(10283,1883)(9083,1883)
\path(2483,83)(2483,683)(3083,1283)
	(4283,83)(2483,83)
\path(4883,83)(5483,1283)(6683,83)(4883,83)
\path(7283,83)(7883,1283)(8483,683)
	(8483,83)(7283,83)
\path(9083,683)(9683,1283)(10283,683)
	(9683,83)(9083,683)
\put(233,308){\makebox(0,0)[lb]{\smash{{\SetFigFont{12}{14.4}{\rmdefault}{\mddefault}{\updefault}$3 \Delta$}}}}
\end{picture}
}
  \caption{All polygons with $i=1$.}
  \label{fig:reflexive}
\end{figure}

Perhaps surprisingly, for $i>0$ no such families exist. For $i=1$,
there are precisely the $16$ lattice equivalence classes depicted in
Figure~\ref{fig:reflexive}. We see that all values $3 \le b \le 9 =
2i+7$ occur. The polygon labeled $3\Delta$ is the $3$-fold dilation of
the standard triangle $\Delta$ which is the convex hull
$\conv[(0,0)\transpose, (1,0)\transpose, (0,1)\transpose]$ of the
origin together with the standard unit vectors. It will play an
important r\^ole later on.

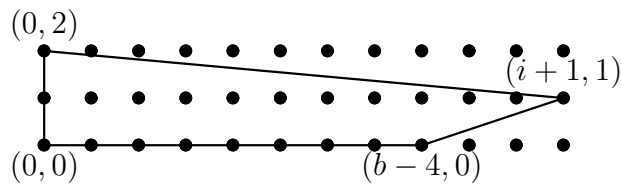
\begin{figure}[htbp]
  \centering
\setlength{\unitlength}{0.00041667in}
\begingroup\makeatletter\ifx\SetFigFont\undefined%
\gdef\SetFigFont#1#2#3#4#5{%
  \reset@font\fontsize{#1}{#2pt}%
  \fontfamily{#3}\fontseries{#4}\fontshape{#5}%
  \selectfont}%
\fi\endgroup%
{\renewcommand{\dashlinestretch}{30}
\begin{picture}(6766,2211)(0,-10)
\put(1883,1083){\blacken\ellipse{150}{150}}
\put(1883,1083){\ellipse{150}{150}}
\put(2483,1083){\blacken\ellipse{150}{150}}
\put(2483,1083){\ellipse{150}{150}}
\put(3083,1083){\blacken\ellipse{150}{150}}
\put(3083,1083){\ellipse{150}{150}}
\put(683,1083){\blacken\ellipse{150}{150}}
\put(683,1083){\ellipse{150}{150}}
\put(1283,1083){\blacken\ellipse{150}{150}}
\put(1283,1083){\ellipse{150}{150}}
\put(3683,1083){\blacken\ellipse{150}{150}}
\put(3683,1083){\ellipse{150}{150}}
\put(83,1083){\blacken\ellipse{150}{150}}
\put(83,1083){\ellipse{150}{150}}
\put(1883,1683){\blacken\ellipse{150}{150}}
\put(1883,1683){\ellipse{150}{150}}
\put(2483,1683){\blacken\ellipse{150}{150}}
\put(2483,1683){\ellipse{150}{150}}
\put(3083,1683){\blacken\ellipse{150}{150}}
\put(3083,1683){\ellipse{150}{150}}
\put(683,1683){\blacken\ellipse{150}{150}}
\put(683,1683){\ellipse{150}{150}}
\put(1283,1683){\blacken\ellipse{150}{150}}
\put(1283,1683){\ellipse{150}{150}}
\put(3683,1683){\blacken\ellipse{150}{150}}
\put(3683,1683){\ellipse{150}{150}}
\put(83,1683){\blacken\ellipse{150}{150}}
\put(83,1683){\ellipse{150}{150}}
\put(1883,483){\blacken\ellipse{150}{150}}
\put(1883,483){\ellipse{150}{150}}
\put(2483,483){\blacken\ellipse{150}{150}}
\put(2483,483){\ellipse{150}{150}}
\put(3083,483){\blacken\ellipse{150}{150}}
\put(3083,483){\ellipse{150}{150}}
\put(683,483){\blacken\ellipse{150}{150}}
\put(683,483){\ellipse{150}{150}}
\put(1283,483){\blacken\ellipse{150}{150}}
\put(1283,483){\ellipse{150}{150}}
\put(3683,483){\blacken\ellipse{150}{150}}
\put(3683,483){\ellipse{150}{150}}
\put(83,483){\blacken\ellipse{150}{150}}
\put(83,483){\ellipse{150}{150}}
\put(6083,1083){\blacken\ellipse{150}{150}}
\put(6083,1083){\ellipse{150}{150}}
\put(6683,1083){\blacken\ellipse{150}{150}}
\put(6683,1083){\ellipse{150}{150}}
\put(4883,1083){\blacken\ellipse{150}{150}}
\put(4883,1083){\ellipse{150}{150}}
\put(5483,1083){\blacken\ellipse{150}{150}}
\put(5483,1083){\ellipse{150}{150}}
\put(4283,1083){\blacken\ellipse{150}{150}}
\put(4283,1083){\ellipse{150}{150}}
\put(6083,1683){\blacken\ellipse{150}{150}}
\put(6083,1683){\ellipse{150}{150}}
\put(6683,1683){\blacken\ellipse{150}{150}}
\put(6683,1683){\ellipse{150}{150}}
\put(4883,1683){\blacken\ellipse{150}{150}}
\put(4883,1683){\ellipse{150}{150}}
\put(5483,1683){\blacken\ellipse{150}{150}}
\put(5483,1683){\ellipse{150}{150}}
\put(4283,1683){\blacken\ellipse{150}{150}}
\put(4283,1683){\ellipse{150}{150}}
\put(6083,483){\blacken\ellipse{150}{150}}
\put(6083,483){\ellipse{150}{150}}
\put(6683,483){\blacken\ellipse{150}{150}}
\put(6683,483){\ellipse{150}{150}}
\put(4883,483){\blacken\ellipse{150}{150}}
\put(4883,483){\ellipse{150}{150}}
\put(5483,483){\blacken\ellipse{150}{150}}
\put(5483,483){\ellipse{150}{150}}
\put(4283,483){\blacken\ellipse{150}{150}}
\put(4283,483){\ellipse{150}{150}}
\thicklines
\path(83,483)(4883,483)(6683,1083)
	(83,1683)(83,483)
\put(83,108){\makebox(0,0)[b]{\smash{{{\SetFigFont{12}{14.4}{\rmdefault}{\mddefault}{\updefault}$(0,0)$}}}}}
\put(6683,1308){\makebox(0,0)[b]{\smash{{{\SetFigFont{12}{14.4}{\rmdefault}{\mddefault}{\updefault}$(i+1,1)$}}}}}
\put(83,1908){\makebox(0,0)[b]{\smash{{{\SetFigFont{12}{14.4}{\rmdefault}{\mddefault}{\updefault}$(0,2)$}}}}}
\put(4883,108){\makebox(0,0)[b]{\smash{{{\SetFigFont{12}{14.4}{\rmdefault}{\mddefault}{\updefault}$(b-4,0)$}}}}}
\end{picture}
}
  \caption{$4\le b\le 2i+6$}
  \label{fig:shoe}
\end{figure}

What can we do for $i \ge 2$? The family shown in
Figure~\ref{fig:shoe} yields all $4\le b\le 2i+6$. In fact, 
Scott~\cite{Scott} showed that $2i+6$ is how far we can get.

Scott's proof is elementary and short enough to be included in this
paper. We give two other proofs for the same result.
One of them uses toric geometry;
it is merely the observation that a well-known
inequality~\cite[Theorem~6]{Schicho:98d} in algebraic geometry
translates to Scott's inequality when applied to toric surfaces.
The third proof is again elementary, and it was the search for this
proof which sparked polygon addiction in the first author.


For the inequality $b\le 2\ i +6$, we have arbitrary large examples
where equality holds (see Figure~\ref{fig:shoe}); but for all these
examples, all interior points are collinear. Under the additional
assumption that the interior points are not collinear, the inequality
can be strengthened to $b \le i + 9$ (see the remark after
Lemma~\ref{lemma:step}).
The coefficient in front of the $i$ can be improved further
by introducing the {\em level} of a 
polygon:
roughly speaking, this is the number of times one can pass
to the convex hull of the interior lattice points.

\begin{figure}[htbp]
  \centering
  \includegraphics{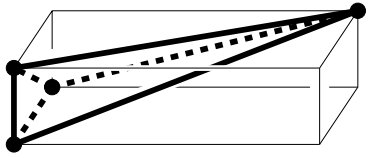}
  \caption{Reeve's simplices}
  \label{fig:reeve}
\end{figure}

Before we really get going, here is a little caveat. Most of our
considerations break down in dimension $3$. Pick's formula has no
analogue. Already tetrahedra with no boundary or interior lattice
points except the vertices can have arbitrary volume. This was first
pointed out by J. Reeve~\cite{Reeve} (see Figure~\ref{fig:reeve}).
Nevertheless, the phenomenon that for given $i>0$, the volume is
bounded occurs in arbitrary dimension~\cite{lagariasZiegler}.

\section{Three proofs of $\mathbf{b \le 2i+7}$}\noindent
Let $P$ be a 
polygon with interior lattice
points. Denote $a$ its surface area, $i$ the number of interior
lattice points, and $b$ the number of lattice points on $P$'s
boundary. In view of Pick's Theorem~\eqref{eq:pick}, the following
three inequalities are equivalent.
\begin{proposition} \label{prop:main}
  If $i>0$, then
%
  \begin{align}
    b \ &\le \ 2 \ i \ + \ 7 \label{eq:main:bi} \\
    a \ &\le \ 2 \ i \ + \ 5/2 \label{eq:main:ai}
    \\
    b \ &\le \ a \ + \ 9/2 \label{eq:main:ab} 
  \end{align}
  with equality only for the triangle $3\Delta$ in
  Figure~\ref{fig:reflexive}.
\end{proposition}
\subsection{Scott's proof} \label{ssec:scott}
\begin{figure}[htbp]
  \centering
\begin{picture}(0,0)%
\includegraphics{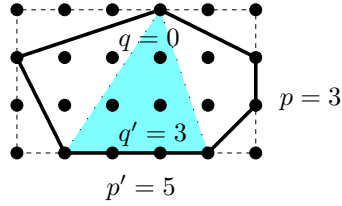}%
\end{picture}%
\setlength{\unitlength}{1973sp}%
\begingroup\makeatletter\ifx\SetFigFont\undefined%
\gdef\SetFigFont#1#2#3#4#5{%
  \reset@font\fontsize{#1}{#2pt}%
  \fontfamily{#3}\fontseries{#4}\fontshape{#5}%
  \selectfont}%
\fi\endgroup%
\begin{picture}(3383,2489)(518,-6967)
\put(3901,-5761){\makebox(0,0)[lb]{\smash{\SetFigFont{10}{12.0}{\rmdefault}{\mddefault}{\updefault}{\color[rgb]{0,0,0}$p=3$}%
}}}
\put(1876,-6211){\makebox(0,0)[lb]{\smash{\SetFigFont{10}{12.0}{\rmdefault}{\mddefault}{\updefault}{\color[rgb]{0,0,0}$q'=3$}%
}}}
\put(1726,-6886){\makebox(0,0)[lb]{\smash{\SetFigFont{10}{12.0}{\rmdefault}{\mddefault}{\updefault}{\color[rgb]{0,0,0}$p'=5$}%
}}}
\put(1876,-5011){\makebox(0,0)[lb]{\smash{\SetFigFont{10}{12.0}{\rmdefault}{\mddefault}{\updefault}{\color[rgb]{0,0,0}$q=0$}%
}}}
\end{picture}
  \caption{$P$ in a box.}
  \label{fig:box}
\end{figure}

Apply lattice equivalences to $P$
so that $P$ fits tightly into a box $[0,p'] \times [0,p]$
with $p$ as small as possible.
Then $2 \le p \le p'$ (remember,
$i>0$). If $P$ intersects the top and the bottom edge of the box in
segments of length $q \ge 0$ and $q' \ge 0$ respectively, then
(See Figure~\ref{fig:box}.)
\begin{align}
  b &\le q+q'+2p \ , \ \text{ and } \label{eq:sc1} \\
  a &\ge p\,(q+q')/2 \ . \label{eq:sc2}
\end{align}

We distinguish three cases
\begin{itemize}
\item[($i$)] $p=2$, or $q+q' \ge 4$, or $p=q+q'=3$
\item[($ii$)] $p=3$, and $q+q' \le 2$
\item[($iii$)] $p \ge 4$, and $q+q' \le 3$.
\end{itemize}
The above inequalities~\eqref{eq:sc1}, and \eqref{eq:sc2} are already
sufficient to deal with the first two cases.
\\[.3\baselineskip]
($i$) We have
\begin{align*}
  2b-2a \ &\le \ 2(q+q'+2p) - p\,(q+q') \\
  &= \ (q+q'-4)(2-p) + 8 \quad \le \ 9 ,
\end{align*}
which shows~\eqref{eq:main:ab} in Proposition~\ref{prop:main}.
(With equality if and only if $p=q+q'=3$, $a=9/2$, $b=9$.)%
\footnote{It is an exercise to show that the only $P$ with these
  parameters is the triangle $3\Delta$ in Figure~\ref{fig:reflexive}.}
\\[.3\baselineskip]
($ii$) The estimate $b \le q+q'+2p \le 8$ together with $i \ge 1$ show
  that inequality~\eqref{eq:main:bi} in Proposition~\ref{prop:main} is
  strictly satisfied.
\\[.3\baselineskip]
($iii$)
The only case where we have to work a little is case three.
Choose points $\pt{x}=(x_1,p)^t$, $\pt{x}'=(x'_1,0)^t$,
$\pt{y}=(0,y_2)^t$, and $\pt{y}'=(p',y'_2)^t$ in $P$ so that $\delta =
|x_1-x'_1|$ is as small as possible. Then $a \ge p (p'- \delta)/2$
(see Figure~\ref{fig:case3}).

\begin{figure}[htbp]
  \centering
\begin{picture}(0,0)%
\includegraphics{case3a.pstex}%
\end{picture}%
\setlength{\unitlength}{1184sp}%
\begingroup\makeatletter\ifx\SetFigFontNFSS\undefined%
\gdef\SetFigFontNFSS#1#2#3#4#5{%
  \reset@font\fontsize{#1}{#2pt}%
  \fontfamily{#3}\fontseries{#4}\fontshape{#5}%
  \selectfont}%
\fi\endgroup%
\begin{picture}(7080,5691)(436,-5329)
\put(7501,-2911){\makebox(0,0)[lb]{\smash{{\SetFigFontNFSS{10}{12.0}{\rmdefault}{\mddefault}{\updefault}{\color[rgb]{0,0,0}$\pt{y}'$}%
}}}}
\put(451,-1111){\makebox(0,0)[lb]{\smash{{\SetFigFontNFSS{10}{12.0}{\rmdefault}{\mddefault}{\updefault}{\color[rgb]{0,0,0}$\pt{y}$}%
}}}}
\put(2776,-61){\makebox(0,0)[lb]{\smash{{\SetFigFontNFSS{10}{12.0}{\rmdefault}{\mddefault}{\updefault}{\color[rgb]{0,0,0}$\pt{x}$}%
}}}}
\put(4051,-5161){\makebox(0,0)[lb]{\smash{{\SetFigFontNFSS{10}{12.0}{\rmdefault}{\mddefault}{\updefault}{\color[rgb]{0,0,0}$\pt{x}'$}%
}}}}
\put(3601,-2536){\makebox(0,0)[b]{\smash{{\SetFigFontNFSS{9}{10.8}{\rmdefault}{\mddefault}{\updefault}{\color[rgb]{0,0,0}\makebox[.3in]{$\overset{\displaystyle\delta}{\longleftrightarrow}$}}%
}}}}
\end{picture}%
\qquad
\begin{picture}(0,0)%
\includegraphics{case3b.pstex}%
\end{picture}%
\setlength{\unitlength}{1184sp}%
\begingroup\makeatletter\ifx\SetFigFont\undefined%
\gdef\SetFigFont#1#2#3#4#5{%
  \reset@font\fontsize{#1}{#2pt}%
  \fontfamily{#3}\fontseries{#4}\fontshape{#5}%
  \selectfont}%
\fi\endgroup%
\begin{picture}(6024,4911)(1189,-5260)
\put(4051,-5161){\makebox(0,0)[lb]{\smash{\SetFigFont{7}{8.4}{\rmdefault}{\mddefault}{\updefault}{\color[rgb]{1,1,1}\mbox{}}%
}}}
\end{picture}
  \caption{\mbox{Case three. Two triangles of total area $p (p'-
      \delta)/2$.}}
  \label{fig:case3}
\end{figure}

Now the task is to apply lattice equivalences so that $\delta$ becomes
small.
\begin{exercise}
  After applying a lattice equivalence of the form $\left[
    \begin{smallmatrix}
      1&k\\0&1
    \end{smallmatrix}
  \right]$
  it is possible to choose $\delta \le (p-q-q')/2$.
\end{exercise}
This lattice equivalence will leave $q$, $q'$, $p$ unchanged, because
it fixes the $x_1$-axis. We still have $p \le p'$ because $p$ was
supposed to be minimal. Thus, we obtain
\begin{equation} \label{eq:sc3}
  a \ge p(p+q+q')/4 ,
\end{equation}
and
\begin{multline*}
  4(b-a) \le 8p+4q+4q' \ - \ p(p+q+q') \\
  = p(8-p) - (p-4)(q+q') \le p(8-p) \le 16
\end{multline*}
because $p\ge4$ in case three. This proves that
inequality~\eqref{eq:main:ab} in Proposition~\ref{prop:main} is
strictly satisfied.
\qed
\subsection{Clipping off vertices}
This proof proceeds by induction on $i$. If $i=1$, we can check the
inequalities on all $16$ lattice equivalence classes of such
$P$. (See Figure~\ref{fig:reflexive}.) 

For the induction step, we want to ``chop off a vertex''.
If $i \ge 2$, and $b \le 10$, nothing is to show.
So assume $b \ge 11$.
By applying a lattice equivalence, we may assume without loss of
generality that $\nv$ and $(1,0)\transpose$ lie in the interior of
$P$. Reflect in the $x_1$-axis if necessary in order to assure
that there are $\ge 5$ boundary lattice points with positive second
coordinate.

First, suppose there is a vertex $v$ with positive second coordinate
which is not unimodular. That is, the triangle formed by $v$ together
with its two neighboring lattice points $v'$ and $v''$ on the boundary
has area $> 1/2$. Denote $P'$ the convex hull $\conv(P \cap \Z^2
\setminus \{v\})$. This omission affects our parameters as follows:
$b' = b+k-2$, $i' = i - k + 1$, and, by Pick's formula, $a' = a-k/2$.
Here $k$ is the lattice length of the boundary of $P'$ that is visible
from $v$. Because $v$ was not unimodular, there is an additional
lattice point in the triangle $vv'v''$. Thus, we have $k \ge
2$. Because there are other lattice points with positive second
coordinate, at least one of $\nv$ or $(1,0)\transpose$ remains in the
interior of $P'$, and we can use induction.

\begin{figure}[htbp]
  \centering
  \includegraphics{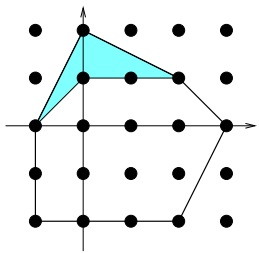}
  \caption{Clipping off a nonunimodular vertex.}
  \label{fig:omit1}
\end{figure}

Now, if all vertices with positive second coordinate are unimodular,
similiarly omit one vertex $v$ together with its two boundary
neighbors $v'$ and $v''$: $P'=\conv(P \cap \Z^2 \setminus
\{v,v',v''\})$.
The parameters change as follows:
$b' = b+k-4$, $i' = i - k + 1$, and $a' = a-k/2-1$, where $k$ is
the lattice length of the boundary of $P'$ that is visible from the
removed points.
In order to see that $k \ge 2$, observe that the point $v'''=v'+v''-v$
belongs to the interior of $P$, and the two adjacent segments of $P'$
are both visible from the removed points.
As observed above, there remain lattice points with
positive second coordinate in $P'$ so that at least one of
$\nv$ or $(1,0)\transpose$ stays in the interior of $P'$.
\hfill $\Box$

\begin{figure}[htbp]
  \centering
\begin{picture}(0,0)%
\includegraphics{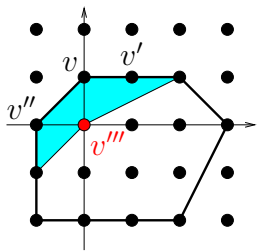}%
\end{picture}%
\setlength{\unitlength}{1973sp}%
\begingroup\makeatletter\ifx\SetFigFont\undefined%
\gdef\SetFigFont#1#2#3#4#5{%
  \reset@font\fontsize{#1}{#2pt}%
  \fontfamily{#3}\fontseries{#4}\fontshape{#5}%
  \selectfont}%
\fi\endgroup%
\begin{picture}(3174,3099)(814,-4348)
\put(2401,-2011){\makebox(0,0)[b]{\smash{{\SetFigFont{12}{14.4}{\rmdefault}{\mddefault}{\updefault}{\color[rgb]{0,0,0}$v'$}%
}}}}
\put(1201,-2686){\makebox(0,0)[rb]{\smash{{\SetFigFont{12}{14.4}{\rmdefault}{\mddefault}{\updefault}{\color[rgb]{0,0,0}$v''$}%
}}}}
\put(1726,-2086){\makebox(0,0)[rb]{\smash{{\SetFigFont{12}{14.4}{\rmdefault}{\mddefault}{\updefault}{\color[rgb]{0,0,0}$v$}%
}}}}
\put(1876,-3136){\makebox(0,0)[lb]{\smash{{\SetFigFont{12}{14.4}{\rmdefault}{\mddefault}{\updefault}{\color[rgb]{1,0,0}$v'''$}%
}}}}
\end{picture}%
  \caption{Clipping off a unimodular vertex (and its neighbors).}
  \label{fig:omit2}
\end{figure}

\subsection{Algebraic geometry} \noindent
We use the letters $d$ and $p$ to denote the
degree and the sectional genus of an algebraic surface. The inequality
$p\le (d-1)(d-2)/2$ holds for arbitrary algebraic surfaces.
If the surface is rational, i.e. if it has a parametrization by
rational functions, then there are more inequalities.

\begin{theorem} \label{thm:ag} \mbox{} 
  \begin{itemize}
  \item If $p=1$, then $d\le 9$.
  \item If $p\ge 2$, then $d\le 4p+4$.
  \end{itemize}
\end{theorem}

Rational surfaces with $p=1$ are called Del Pezzo surfaces. The degree
bound $9$ is due to del Pezzo~\cite{Pezzo:87}. The bound $d\le 4p+4$
was shown by Jung~\cite{Jung:90}, hence this proof is
actually the oldest one! 
A modern proof can be found in Schicho~\cite{Schicho:98d}.

Toric surfaces are rational, and Scott's inequality is equivalent
to Theorem~\ref{thm:ag} for toric surfaces.
\hfill $\Box$

\section{Onion skins} \label{sec:onions}
\begin{figure}[htbp]
  \centering
  \includegraphics[height=43mm]{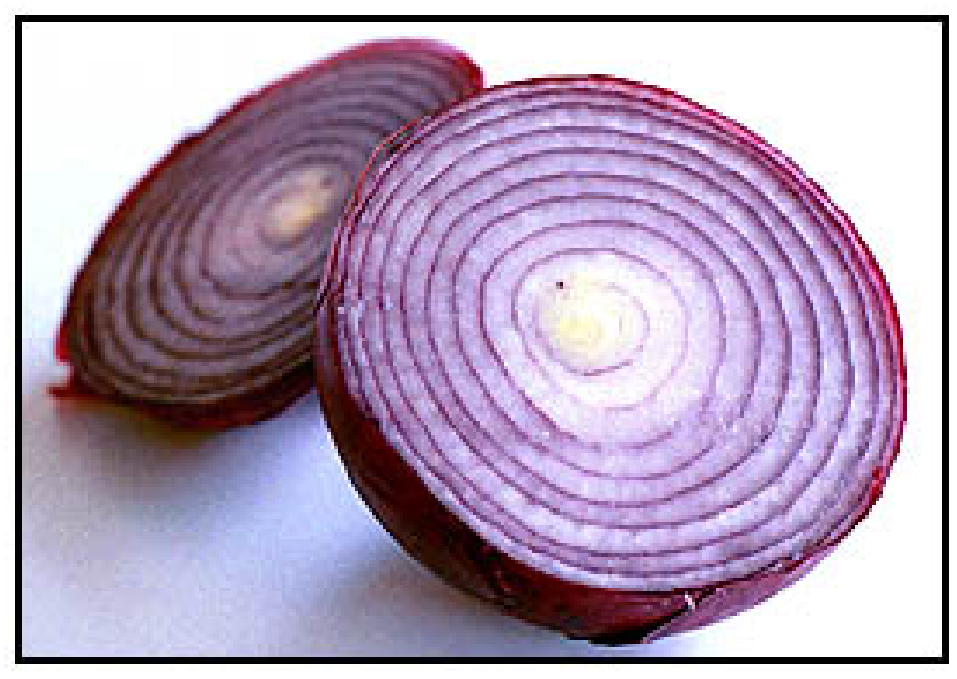}
\end{figure}

In flatland, take a solid polygon $P$ into your hand
and peel off the shell. You get another convex polygon $P^{(1)}$, the
convex hull of its interior lattice points. Except, of course, if $i=0$
then $P$ was an empty nut, and if all interior lattice points
are collinear then $P^{(1)}$ is a ``degenerate polygon'', namely a line
segment or a single point.

Repeat this process as long as it is possible, peeling off
the skins of the polygon one after the other: $P^{(k+1)} :=
(P^{(k)})^{(1)}$. After $n$ steps you arrive at the nucleus $P^{(n)}$,
which is either a degenerate polygon or an empty nut. We define the
{\em level\/} $\ell=\ell(P)$ in the following way:
\begin{itemize}
  \item $\ell(P)=n$ if the nucleus is a degenerate polygon,
  \item $\ell(P)=n+1/3$ if the nucleus is $\Delta$,
  \item $\ell(P)=n+2/3$ if the nucleus is $2\Delta$, and
  \item $\ell(P)=n+1/2$ if the nucleus is any other empty nut.
\end{itemize}
Here $\Delta$ stands for (a polygon lattice equivalent to) the
standard triangle $\conv[(0,0)\transpose, (1,0)\transpose,
(0,1)\transpose]$. The purpose of this weird definition is
to ensure the second statement in the exercise below.

\begin{exercise}
  Show that $\ell$ is uniquely defined by
  \begin{itemize}
  \item $\ell(P) = \ell(P^{(1)})+1$ if $P^{(1)}$ is $2$-dimensional,
    and
  \item $\ell(kP)=k\ell(P)$ for positive integers $k$.
  \end{itemize}


 \end{exercise}

\begin{figure}[htbp]
  \centering
  \includegraphics{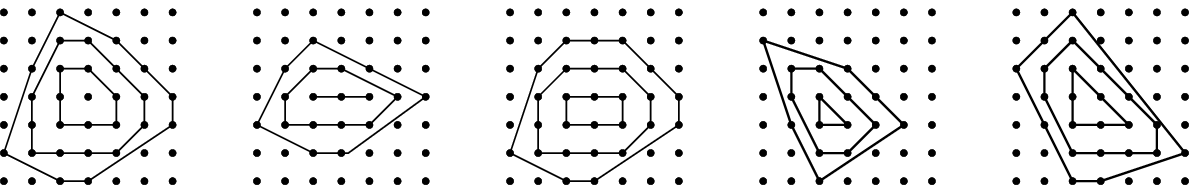}
  \caption{Polygons of levels \hfill \mbox{} \mbox{$\ell=3$, $\ell=2$,
      $\ell=5/2$, $\ell=7/3$, and $\ell=8/3$.}}
  \label{fig:P(ell)}
\end{figure}

The level of a polygon is an analogue of the radius of the in-circle
in Euclidean geometry. There we have the equation $2a=\ell b$.
In lattice geometry, we have an inequality.

\begin{oniontheorem} \stepcounter{theorem} \label{thm:onion}
  Let $P$ be a convex lattice polygon of area $a$ and
  level $\ell \ge 1$ with $b$
  and $i$ lattice points on the boundary and in the interior,
  respectively. Then $(2\ell-1) b \le 2i + 9\ell^2 - 2$, or
  equivalently $2\ell b\le 2a+9\ell^2$, or equivalently
  $(4l-2)a\le 9\ell^2+4l(i-1)$, with equality
  if and only if $P$ is a multiple of $\Delta$.
\end{oniontheorem}
For $\ell>1$, these inequalities really strengthen the old $b \le
2i+7$.
We give two elementary proofs. One is similar to Scott's proof.
The other is a bit longer, but it gives more insight into the process
of peeling onion skins. For instance, it reveals that the set of 
all polygons $P$ such that $P^{(1)}=Q$ for some fixed $Q$ is either
empty or has a largest element.
\subsection{Moving out edges}
Using this technique, it is actually possible to sharpen the bound 
in various (sub)cases. E.g.,
\begin{itemize}
  \item if $P^{(\ell)}$ is a point, but $P^{(\ell-1)} \neq 3 \Delta$,
  then $(2\ell-1) b \le 2i + 8 \ell^2 - 2$;
  \item if $P^{(\ell)}$ is a segment, then $(2\ell-1) b \le 2i +
  8\ell^2 - 2$ with equality if and only if $P$ is lattice equivalent
  to a polygon with vertices $\nv, (r,0)\transpose,
  (2pq+r,2p)\transpose, (0,2p)\transpose$ for integers $p \ge 1$, $q,r
  \ge 0$ such that $pq+r \ge 3$;
  \item if $P^{(\ell)}$ has no interior lattice points but is not a
  multiple of $\Delta$, then $(2\ell-1) b \le 2i + 8 \ell^2 - 2$.
\end{itemize}
We reduce the proof to the case that $P$ is obtained from $P^{(1)}$ by
``moving out the edges by one''. This is done in the following three
lemmas. Finally, Lemma~\ref{lemma:step} yields the induction step in
the proof of the Onion--Skin Theorem.

We say that an inequality $\langle \pt{a}, \pt{x} \rangle = a_1 x_1 +
a_2 x_2 \le b$ with coprime $(a_1,a_2)$ defines an edge of a 
polygon $Q$ if it is satisfied by all points $\pt{x} \in Q$, and there
are two distinct points in $Q$ satisfying equality.
Then, moving out this edge by one means to relax the inequality to
$\langle \pt{a}, \pt{x} \rangle \le b + 1$.
\begin{lemma} \label{lemma:shift}
  Suppose that the inequality $\langle \pt{a}, \pt{x} \rangle \le b$
  defines an edge of $P^{(1)}$. Then $\langle \pt{a}, \pt{x} \rangle
  \le b + 1$ is valid for $P$.
\end{lemma}
That means, if we move all the edges of $Q=P^{(1)}$ out by one, we
obtain a superset $Q^{(-1)}$ of $P$.

\begin{figure}[htbp]
  \centering
  \includegraphics{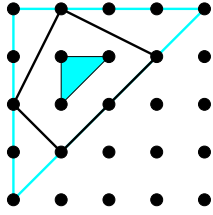}
  \caption{If $Q=P^{(1)}$ then $P \subseteq Q^{(-1)}$.}
  \label{fig:shift}
\end{figure}

\begin{proof}
  We may apply a lattice equivalence to reduce to the case
  where the edge is defined by $x_2 \le 0$, and that $(0,0)\transpose$
  and $(1,0)\transpose$ are two lattice points of $P^{(1)}$ lying
  on this edge. 
  Suppose indirectly that 
  $P$ has a vertex $\pt{v}$ with $v_2 > 1$. Then the triangle formed
  by the three points has area $v_2/2 \ge 1$. It must therefore
  contain another lattice point which lies in the interior of $P$, and
  has positive second coordinate.
\end{proof}
For arbitrary $Q$, this $Q^{(-1)}$ does not necessarily have integral
vertices.
But then, not every 
polygon arises as $P^{(1)}$ for some $P$. A necessary condition is
that the polygon has good angles.

\begin{figure}[htbp]
  \centering
  \includegraphics{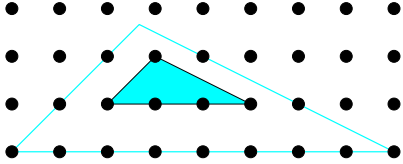}
  \caption{$Q^{(-1)}$ may be nonintegral.}
  \label{fig:badshift}
\end{figure}

\begin{lemma} \label{lemma:angles}
  If $P^{(1)}$ is $2$--dimensional, then for all vertices $\pt{v}$ of
  $P^{(1)}$, the cones generated by $P^{(1)} - \pt{v}$ are
  lattice equivalent to a cone generated by $(1,0)\transpose$ and
  $(-1,k)\transpose$, for some integer $k \ge 1$.
\end{lemma}

\begin{figure}[htbp]
  \centering
  \includegraphics{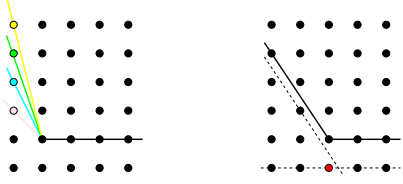}
  \caption{Good angles and a bad angle.}
  \label{fig:goodAngles}
\end{figure}

\begin{proof}
  Assume, after a lattice equivalence, that $\pt{v} = \nv$, and the
  rays of the cone in question are generated by $(1,0)\transpose$ and
  $(-p,q)\transpose$, with coprime $0 < p \le q$ (see
  Exercise~\ref{ex:angle}). By Lemma~\ref{lemma:shift}, all points of
  $P$ satisfy $x_2 \ge -1$ and $q x_1 + p x_2 \ge -1$. But this
  implies $x_1 + x_2 \ge -1 + \frac{p-1}{q}$. So, if $p > 1$, because
  $P$ has integral vertices, we have $x_1 + x_2 \ge 0$ for all points
  of $P$. This contradicts the fact that $\nv \in P^{(1)}$.
\end{proof}
For a vertex $\pt{v}$ of a 
polygon define the shifted
vertex $\pt{v}^{(-1)}$ as follows. Let $\langle \pt{a}, \pt{x} \rangle
\le b$ and \mbox{$\langle \pt{a}', \pt{x} \rangle \le b'$} be the two
edges that intersect in $\pt{v}$.
The unique solution to $\langle \pt{a}, \pt{x} \rangle = b
+ 1$ and $\langle \pt{a}', \pt{x} \rangle = b' + 1$ is denoted
$\pt{v}^{(-1)}$. According to Lemma~\ref{lemma:angles}, when we deal
with $P^{(1)}$ then $\pt{v}^{(-1)}$ is a lattice point. (In the
situation of the lemma, it is $(0,-1)\transpose$.)
We obtain a characterization of when $Q = P^{(1)}$ for some $P$.
\begin{lemma} \label{lemma:Qminus1integral}
  For a 
polygon $Q$, the following are equivalent:
  \begin{itemize}
    \item $Q = P^{(1)}$ for some 
      polygon $P$.
    \item $Q^{(-1)}$ has integral vertices.
  \end{itemize}
\end{lemma}
Thus, given $Q$, {\em the\/} maximal polygon $P$ with $P^{(1)} = Q$ is
$P = Q^{(-1)}$. We will (and can) restrict to this situation when we
prove the induction step $\ell \leadsto \ell + 1$ for the Onion--Skin
Theorem.
\begin{proof}
  If $Q^{(-1)}$ has integral vertices, then its interior lattice
  points span $Q$.
  For the converse direction, if $Q = P^{(1)}$ then we claim that
  \begin{equation} \label{eq:Qminus1}
    Q^{(-1)} = \conv \{ \pt{v}^{(-1)} \suchthat \pt{v} \text{ vertex
      of } Q \} \ .
  \end{equation}
  To this end,
  denote $\pt{a_1}, \ldots, \pt{a_n}$ and $b_1, \ldots, b_n$ the
  normal vectors respectively right hand sides of the edges of $Q$ in
  cyclic order. Also, denote $\pt{v}_1, \ldots, \pt{v}_n$ the vertices
  of $Q$ so that edge number $k$ is the segment
  $\segment{v_\text{$k$}}{v_\text{$k\negthinspace+\negthinspace1$}}$
  ($k \mod n$).
  \smallskip

  \noindent ``$\subseteq$'': this inclusion holds for
  arbitrary $Q$\/.
  For a point $\pt{y} \in Q^{(-1)}$, let $\langle \pt{a}_k, \pt{x}
  \rangle \le b_k$ be an edge of $Q$ that maximizes $\langle \pt{a},
  \pt{y} \rangle - b$ over all edges. So if $\langle \pt{a}_k, \pt{y}
  \rangle - b_k \le 0$ then $\pt{y} \in Q$. Otherwise we have
  \begin{itemize}
  \item $b_k \le \langle \pt{a}_k, \pt{y} \rangle \le b_k + 1$,
  \item $\langle \pt{a}_k, \pt{y} \rangle - b_k \ge \langle
    \pt{a}_{k-1}, \pt{y} \rangle - b_{k-1}$, and
  \item $\langle \pt{a}_k, \pt{y} \rangle - b_k \ge \langle
    \pt{a}_{k+1}, \pt{y} \rangle - b_{k+1}$,
  \end{itemize}
  which describes (a subset of) the convex hull of $\pt{v}_k,
  \pt{v}_k^{(-1)}, \pt{v}_{k+1}, \pt{v}_{k+1}^{(-1)}$.
  \smallskip

  \noindent ``$\supseteq$'': For this inclusion we actually use $Q =
  P^{(1)}$. Figure~\ref{fig:Qminus1} shows how
  Equation~\eqref{eq:Qminus1} can fail in general.
  
  \begin{figure}[htb]
    \centering
\begin{picture}(0,0)%
\includegraphics{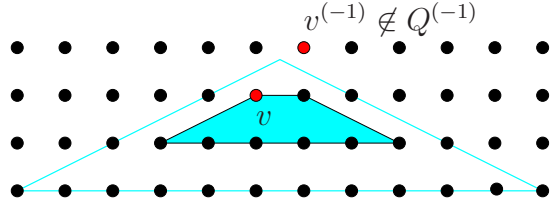}%
\end{picture}%
\setlength{\unitlength}{1973sp}%
\begingroup\makeatletter\ifx\SetFigFont\undefined%
\gdef\SetFigFont#1#2#3#4#5{%
  \reset@font\fontsize{#1}{#2pt}%
  \fontfamily{#3}\fontseries{#4}\fontshape{#5}%
  \selectfont}%
\fi\endgroup%
\begin{picture}(6766,2435)(2918,-7044)
\put(6601,-4936){\makebox(0,0)[lb]{\smash{{\SetFigFont{12}{14.4}{\rmdefault}{\mddefault}{\updefault}$v^{(-1)} \not\in Q^{(-1)}$}}}}
\put(6001,-6136){\makebox(0,0)[lb]{\smash{{\SetFigFont{12}{14.4}{\rmdefault}{\mddefault}{\updefault}$v$}}}}
\end{picture}%
    \caption{Equation \eqref{eq:Qminus1} can fail in general.}
    \label{fig:Qminus1}
  \end{figure}
  
  In our situation, $Q = P^{(1)}$, and we need to show that $v_k^{(-1)}$
  satisfies all inequalities $\langle \pt{a}_j, \cdot \rangle \le b_j +
  1$ for $Q^{(-1)}$.
  Our assumption implies that
  $P$ (and therefore $Q^{(-1)}$ by Lemma~\ref{lemma:shift}) contains
  points $\pt{w}_j$ with $\langle \pt{a}_j, \pt{w}_j \rangle = b_j + 1$.
  None of the other edge normals belongs to the cone
  generated by $\pt{a}_{k-1}$ and $\pt{a}_k$. So for $j \neq k,k-1$,
  \begin{equation*}
    \begin{array}{lcccccc}
      \text{either } & \langle \pt{a}_j, \pt{v}_k^{(-1)} \rangle &\le&
      \langle \pt{a}_j, \pt{w}_{k-1} \rangle &\le& b_m + 1 , \\
      \text{or } & \langle \pt{a}_j, \pt{v}_k^{(-1)} \rangle &\le& \langle
      \pt{a}_j, \pt{w}_k \rangle &\le& b_m + 1 & \text{ (or both).}
    \end{array}
  \end{equation*}
\end{proof}
Finally, we can prove the key lemma for our induction step.
\begin{lemma} \label{lemma:step}
  Let $b^{(1)}$ denote the number of lattice points on the boundary of
  $P^{(1)}$. Then $b \le b^{(1)} + 9$, with equality if and only if
  $P$ is a multiple of $\Delta$.
\end{lemma}
This immediately shows that $b \le b^{(1)} + 9 \le i + 9$ if $P^{(1)}$
is $2$--dimensional.

For the proof, we need a result of B.~Poonen and
F.~Rodriguez-Villegas~\cite{Poonen12}. 
Consider a primitive oriented segment $s=\segment{x}{y}$, i.e.,
$\pt{x}$ and $\pt{y}$ are the only lattice points $s$ contains. Call
$s$ admissible if the triangle $\conv(\nv, \pt{x}, \pt{y})$ contains
no other lattice points. Equivalently, $s$ is admissible if the
determinant $\sign(s) = \left| \begin{smallmatrix} x_1&y_1 \\ x_2&y_2
  \end{smallmatrix} \right|$ is equal to $\pm 1$.
The length of a sequence $(s^{(1)}, \ldots, s^{(n)})$ of admissible
segments is $\sum \sign(s^{(k)})$.

The dual of an admissible segment is the unique integral normal vector
$\pt{a}=\pt{a}(s)$ such that $\langle \pt{a}, \pt{x} \rangle = \langle
\pt{a}, \pt{y} \rangle = 1$. For a closed polygon with segments
$(s^{(1)}, \ldots, s^{(n)})$, the dual polygon walks through the
normal vectors $\pt{a}(s^{(k)})$.

\begin{figure}[htbp]
  \begin{center}
\setlength{\unitlength}{0.00041667in}
\begingroup\makeatletter\ifx\SetFigFontNFSS\undefined%
\gdef\SetFigFontNFSS#1#2#3#4#5{%
  \reset@font\fontsize{#1}{#2pt}%
  \fontfamily{#3}\fontseries{#4}\fontshape{#5}%
  \selectfont}%
\fi\endgroup%
{\renewcommand{\dashlinestretch}{30}
\begin{picture}(5570,2581)(0,-10)
\put(4515,1883){\blacken\ellipse{150}{150}}
\put(4515,1883){\ellipse{150}{150}}
\put(4515,683){\blacken\ellipse{150}{150}}
\put(4515,683){\ellipse{150}{150}}
\put(3315,1283){\blacken\ellipse{150}{150}}
\put(3315,1283){\ellipse{150}{150}}
\put(3315,1883){\blacken\ellipse{150}{150}}
\put(3315,1883){\ellipse{150}{150}}
\put(3315,683){\blacken\ellipse{150}{150}}
\put(3315,683){\ellipse{150}{150}}
\put(3915,1283){\blacken\ellipse{150}{150}}
\put(3915,1283){\ellipse{150}{150}}
\put(3915,1883){\blacken\ellipse{150}{150}}
\put(3915,1883){\ellipse{150}{150}}
\put(3915,683){\blacken\ellipse{150}{150}}
\put(3915,683){\ellipse{150}{150}}
\put(3915,2483){\blacken\ellipse{150}{150}}
\put(3915,2483){\ellipse{150}{150}}
\put(4515,2483){\blacken\ellipse{150}{150}}
\put(4515,2483){\ellipse{150}{150}}
\put(3315,2483){\blacken\ellipse{150}{150}}
\put(3315,2483){\ellipse{150}{150}}
\put(5115,2483){\blacken\ellipse{150}{150}}
\put(5115,2483){\ellipse{150}{150}}
\put(5115,1883){\blacken\ellipse{150}{150}}
\put(5115,1883){\ellipse{150}{150}}
\put(5115,1283){\blacken\ellipse{150}{150}}
\put(5115,1283){\ellipse{150}{150}}
\put(5115,683){\blacken\ellipse{150}{150}}
\put(5115,683){\ellipse{150}{150}}
\put(5115,83){\blacken\ellipse{150}{150}}
\put(5115,83){\ellipse{150}{150}}
\put(4515,83){\blacken\ellipse{150}{150}}
\put(4515,83){\ellipse{150}{150}}
\put(3915,83){\blacken\ellipse{150}{150}}
\put(3915,83){\ellipse{150}{150}}
\put(3315,83){\blacken\ellipse{150}{150}}
\put(3315,83){\ellipse{150}{150}}
\put(1515,1283){\blacken\ellipse{150}{150}}
\put(1515,1283){\ellipse{150}{150}}
\put(1515,1883){\blacken\ellipse{150}{150}}
\put(1515,1883){\ellipse{150}{150}}
\put(1515,683){\blacken\ellipse{150}{150}}
\put(1515,683){\ellipse{150}{150}}
\put(315,1283){\blacken\ellipse{150}{150}}
\put(315,1283){\ellipse{150}{150}}
\put(315,1883){\blacken\ellipse{150}{150}}
\put(315,1883){\ellipse{150}{150}}
\put(315,683){\blacken\ellipse{150}{150}}
\put(315,683){\ellipse{150}{150}}
\put(915,1883){\blacken\ellipse{150}{150}}
\put(915,1883){\ellipse{150}{150}}
\put(915,683){\blacken\ellipse{150}{150}}
\put(915,683){\ellipse{150}{150}}
\put(4515,1283){\whiten\ellipse{150}{150}}
\put(4515,1283){\ellipse{150}{150}}
\put(915,1283){\whiten\ellipse{150}{150}}
\put(915,1283){\ellipse{150}{150}}
\thicklines
\path(1515,1283)(315,1883)
\path(556.495,1829.334)(315.000,1883.000)(502.830,1722.003)
\path(315,1883)(915,1883)
\path(675.000,1823.000)(915.000,1883.000)(675.000,1943.000)
\path(915,1883)(315,683)
\path(368.666,924.495)(315.000,683.000)(475.997,870.830)
\path(315,683)(1515,1283)
\path(1327.170,1122.003)(1515.000,1283.000)(1273.505,1229.334)
\dashline{150.000}(5115,2483)(4515,1883)
\path(4642.279,2095.132)(4515.000,1883.000)(4727.132,2010.279)
\dashline{150.000}(4515,1883)(3315,1883)
\path(3555.000,1943.000)(3315.000,1883.000)(3555.000,1823.000)
\dashline{150.000}(3315,1883)(5115,83)
\path(4902.868,210.279)(5115.000,83.000)(4987.721,295.132)
\dashline{150.000}(5115,83)(5115,2483)
\path(5175.000,2243.000)(5115.000,2483.000)(5055.000,2243.000)
\put(1140,1583){\makebox(0,0)[lb]{\smash{{\SetFigFontNFSS{12}{14.4}{\rmdefault}{\mddefault}{\updefault}1}}}}
\put(540,2108){\makebox(0,0)[lb]{\smash{{\SetFigFontNFSS{12}{14.4}{\rmdefault}{\mddefault}{\updefault}2}}}}
\put(15,833){\makebox(0,0)[lb]{\smash{{\SetFigFontNFSS{12}{14.4}{\rmdefault}{\mddefault}{\updefault}3}}}}
\put(540,308){\makebox(0,0)[lb]{\smash{{\SetFigFontNFSS{12}{14.4}{\rmdefault}{\mddefault}{\updefault}4}}}}
\put(2865,1958){\makebox(0,0)[lb]{\smash{{\SetFigFontNFSS{12}{14.4}{\rmdefault}{\mddefault}{\updefault}3}}}}
\put(5340,158){\makebox(0,0)[lb]{\smash{{\SetFigFontNFSS{12}{14.4}{\rmdefault}{\mddefault}{\updefault}4}}}}
\put(4215,1958){\makebox(0,0)[lb]{\smash{{\SetFigFontNFSS{12}{14.4}{\rmdefault}{\mddefault}{\updefault}2}}}}
\put(5340,2258){\makebox(0,0)[lb]{\smash{{\SetFigFontNFSS{12}{14.4}{\rmdefault}{\mddefault}{\updefault}1}}}}
\end{picture}
}
    \caption{A polygon and its dual. \hfill \mbox{}
      \mbox{Their lengths are $1-1+1+1$ respectively $1+2+3+4$.}
      }
    \label{fig:dual}
  \end{center}
\end{figure}
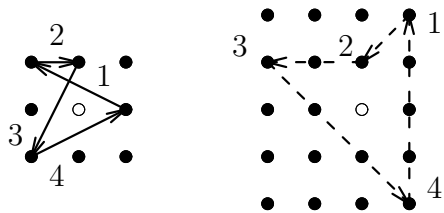

\begin{theorem}[Poonen and Rodriguez-Villegas~\cite{Poonen12}]
  The sum of the lengths of an admissible polygon and its dual is $12$
  times the winding number of the polygon.
\end{theorem}

Heuristically, the winding number counts how many times a polygon
winds around the origin. Dual polygons will have equal winding
number. In this article, we will only be concerned with polygons of
winding number one.

\begin{figure}[htbp]
  \begin{center}
\setlength{\unitlength}{0.00041667in}
\begingroup\makeatletter\ifx\SetFigFont\undefined%
\gdef\SetFigFont#1#2#3#4#5{%
  \reset@font\fontsize{#1}{#2pt}%
  \fontfamily{#3}\fontseries{#4}\fontshape{#5}%
  \selectfont}%
\fi\endgroup%
{\renewcommand{\dashlinestretch}{30}
\begin{picture}(1366,1381)(0,-10)
\put(1283,683){\blacken\ellipse{150}{150}}
\put(1283,683){\ellipse{150}{150}}
\put(1283,1283){\blacken\ellipse{150}{150}}
\put(1283,1283){\ellipse{150}{150}}
\put(1283,83){\blacken\ellipse{150}{150}}
\put(1283,83){\ellipse{150}{150}}
\put(83,683){\blacken\ellipse{150}{150}}
\put(83,683){\ellipse{150}{150}}
\put(83,1283){\blacken\ellipse{150}{150}}
\put(83,1283){\ellipse{150}{150}}
\put(83,83){\blacken\ellipse{150}{150}}
\put(83,83){\ellipse{150}{150}}
\put(683,683){\blacken\ellipse{150}{150}}
\put(683,683){\ellipse{150}{150}}
\put(683,1283){\blacken\ellipse{150}{150}}
\put(683,1283){\ellipse{150}{150}}
\put(683,83){\blacken\ellipse{150}{150}}
\put(683,83){\ellipse{150}{150}}
\thicklines
\path(83,683)(1283,1283)
\path(1095.170,1122.003)(1283.000,1283.000)(1041.505,1229.334)
\path(1283,1283)(683,83)
\path(736.666,324.495)(683.000,83.000)(843.997,270.830)
\path(683,83)(83,1283)
\path(243.997,1095.170)(83.000,1283.000)(136.666,1041.505)
\path(83,1283)(1283,683)
\path(1041.505,736.666)(1283.000,683.000)(1095.170,843.997)
\path(1283,683)(83,83)
\path(270.830,243.997)(83.000,83.000)(324.495,136.666)
\path(83,83)(83,683)
\path(143.000,443.000)(83.000,683.000)(23.000,443.000)
\end{picture}
}
    \caption{A polygon of winding number $-2$.}
    \label{fig:winding}
  \end{center}
\end{figure}
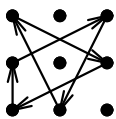

\begin{proof}[Proof of Lemma~\ref{lemma:step}.]
  Let $Q = P^{(1)}$, and note that by Lemma~\ref{lemma:shift}, $P
  \subseteq Q^{(-1)}$, and by Lemma~\ref{lemma:Qminus1integral},
  $Q^{(-1)}$ has integral vertices. Notice that the number of boundary
  lattice points of $Q$ is $b^{(1)}$, and let $b'$ be the number of
  boundary lattice points of $Q^{(-1)}$. Since $P$ and $Q^{(-1)}$ have
  the same interior lattice points and $P \subseteq Q^{(-1)}$, by
  Pick's Theorem $Q^{(-1)}$ has at least as many boundary lattice
  points as $P$; in other words, $b' \ge b$.

  For each of the vertices $v^{(1)}_1, \ldots, v^{(1)}_n$ of $Q$
  there is a corresponding vertex $v_1, \ldots, v_n$ of $Q^{(-1)}$.
  Consider the (possibly nonconvex, nonsimple) admissible polygon
  with vertices
  $v_1-v^{(1)}_1, \ldots, v_n-v^{(1)}_n$. It is admissible because
  there are no lattice points between $Q$ and $Q^{(-1)}$. One can
  think of it as what remains of $Q^{(-1)}$ when $Q$ shrinks to a
  point. Each segment measures the difference (with the correct sign)
  between the corresponding edges of $Q^{(-1)}$ and $Q$. I.e., the
  length of that polygon is precisely $b'-b^{(1)}$.
  
  \begin{figure}[htbp]
    \centering
    \includegraphics{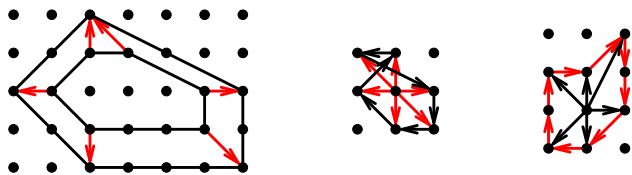}
    \caption{An admissible polygon from $(P,P^{(1)})$, and its dual.}
    \label{fig:shrink}
  \end{figure}

  Now the dual polygon will walk through the normal vectors of
  $Q$. Therefore all segments will count with positive length,
  and there cannot be less than $3$. Also, there is a unique one with
  $3$ segments, which is the dual to $3 \Delta$. Thus $b-b^{(1)} \le
  b'-b^{(1)} \le 12-3$ with equality only for multiples of $\Delta$.
\end{proof}

\begin{proof}[Proof of the Onion--Skin Theorem.]
  Induction on $\ell$. 
  \begin{itemize}
    \item For $\ell=1$, the inequality $b \le 2i + 7$ was proved
    earlier.
    \item For $\ell=4/3$, we have $i=3$, and $P \subseteq 4
    \Delta$. So $b \le 12$.
    \item For $\ell=5/3$, we have $i=6$, and $P \subseteq 5
    \Delta$. So $b \le 15$.
    \item For $\ell=3/2$, Lemma~\ref{lemma:step} reads $b \le i + 8$
    which is stronger than what we need.
  \end{itemize}
  
  If $\ell \ge 2$, we have
  \begin{align*}
    (2\ell-1) b &\le (2\ell-1) b^{(1)} + 9(2\ell-1) \\
    &= 2 b^{(1)} + (2(\ell-1)-1) b^{(1)} + 9(2\ell-1) \\
    &\le 2 b^{(1)} + 2i^{(1)} + 9(\ell-1)^2 - 2 + 9(2\ell-1) \\
    &= 2i + 9\ell^2 - 2
  \end{align*}
\end{proof}

\subsection{Generalizing Scott's proof}
As in Subsection~\ref{ssec:scott}, we tightly fit $P$ into a box
$[0,p'] \times [0,p]$, with $p\le p'$. Let $q$ and $q'$ be the
length of the top and bottom edge (see Figure~\ref{fig:box}).
We again apply lattice equivalence transformations such that $p$ is
as small as possible, and that $P$ has points on the top and bottom
edges with horizontal distance smaller than or equal to $(p-q-q')/2$.
Again, we obtain the following inequalities:
\begin{align*}
b \le q+q'+2p \tag{\ref{eq:sc1}} \\
a \ge p(q+q')/2 \tag{\ref{eq:sc2}} \\
a \ge p(p+q+q')/4 \tag{\ref{eq:sc3}}
\end{align*}
Set $x:=p/\ell$ and $y:=(q+q')/\ell$. Then $x\ge 2$, because
passing to $P^{(1)}$ reduces the height at least by
2.\footnote{We also have $x \le 3$ with equality only for multiples of
  $\Delta$.}
From \eqref{eq:sc1} and \eqref{eq:sc2}, we get
\begin{align*}
\frac{2\ell b-2a-9\ell^2}{\ell^2} \le 2(q+q'+2p)/\ell - p(q+q')/\ell^2-9 \\
        = -xy+4x+2y-9 ,
\end{align*}
and from \eqref{eq:sc1} and \eqref{eq:sc3}, we get
\begin{align*}
\frac{4\ell b-4a-18\ell^2}{\ell^2} \le 4(q+q'+2p)/\ell - p(p+q+q')/\ell^2-18 \\
        = -x^2-xy+8x+4y-18.
\end{align*}

\begin{figure}[htbp]
  \begin{center}
    \includegraphics{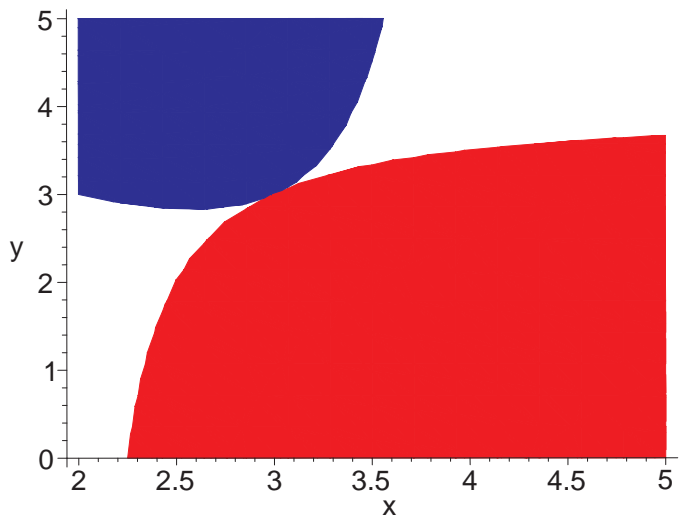}
    \caption{At least one of the polynomials is $\le 0$.}
    \label{fig:plot}
  \end{center}
\end{figure}

For $x\ge 2$ and $y\ge 0$, at least one of the two polynomials
$p_1(x,y) = -xy+4x+2y-9$ and $p_2(x,y) = -x^2-xy+8x+4y-18$ is zero or
negative, as it can be seen in Figure~\ref{fig:plot}. (The two shaded
regions are where $p_1$ and $p_2$, respectively, take non-negative 
values.) There is only
one point where both upper bounds reach zero, namely $(x,y)=(3,3)$,
and this is the only case where equality can hold in the Onion--Skin
Theorem. It is an exercise to show that equality actually holds only 
for multiples of $\Delta$.

\section{Conclusion} \label{sec:conclusion} \noindent
\subsection{Summary of results.}
For a triple $(a,b,i)$ of numbers the following are equivalent.
\begin{itemize}
\item There is a convex lattice polygon $P$ with
  $(a,b,i)=(a(P),b(P),i(P))$.
\item $b \in \Z_{\ge 3}$, $i \in \Z_{\ge 0}, a = i + b/2 -1$, and
  \begin{itemize}
  \item[\scriptsize$\blacktriangleright$] $i=0$ or
  \item[\scriptsize$\blacktriangleright$] $i=1$ and $b \le 9$ or
  \item[\scriptsize$\blacktriangleright$] $i \ge 2$ and $b \le 2i+6$.
  \end{itemize}
\end{itemize}
Furthermore, if $\ell = \ell(P)$, then $(2\ell-1) b \le 2i + 9\ell^2 -
2$. 

\subsection{Outlook.}
Is there a proof of the Onion--Skin Theorem using algebraic geometry?
Currently not (yet). The toric dictionary between 
polygons and
algebraic varieties also does not (yet) have an algebraic geometry term
for the level of a polygon. A first step in this direction is the
use of the {\em process\/} of peeling off onion skins -- or rather
its algebraic geometry analogue -- for the simplification of
the rational parametrization of an algebraic surface~\cite{Schicho:03a}.

In any case, the Onion--Skin Theorem gives rise to a {\em conjecture\/}
in algebraic geometry, namely the inequality 
$(2\ell-1)d\le 9\ell^2+4\ell(p-1)$ for any rational surface
of degree $d$, sectional genus $p$, and level $\ell$. Here the level
of an algebraic surface should be defined via the process of peeling
mentioned above. For toric surfaces, the inequality holds
by the Onion--Skin Theorem, but for nontoric rational surfaces
we do not have a proof (nor a counterexample).

We thank Ricky Pollack for his encouragement to write up this story.
We also thank Daniel J. Velleman and two anonymous referees who helped
us do so in a comprehensible way. 

\bibliographystyle{plain}
\bibliography{alles,josef}
\setlength{\parindent}{0pt}

\end{document}